\documentclass[reqno,a4paper,11pt]{amsart}

\usepackage[all]{xy}
\usepackage{amsfonts}
\usepackage[mathcal]{eucal}
\usepackage{eufrak}
\usepackage{amssymb}
\usepackage{amsmath}
\usepackage{mathrsfs}
\usepackage{color}
\usepackage[pagebackref,colorlinks]{hyperref}
\usepackage[pagebackref]{hyperref}
\usepackage{enumerate}

\newcommand{\ra}{\rightarrow}

\newcommand{\mi}{\setminus}

\newcommand{\ov}{\overline}

\newcommand{\cc}{\mathbb C}
\newcommand{\qq}{\mathbb Q}
\newcommand{\rr}{\mathbb R}
\newcommand{\zz}{\mathbb Z}
\newcommand{\nn}{\mathbb N}

\newcommand{\dlambda}{\partial}
\newcommand{\id}{\operatorname{id}}

\newcommand{\sE}{\mathsf{E}}
\newcommand{\sT}{\mathsf{T}}
\newcommand{\sR}{\mathsf{R}}

\newcommand{\sL}{\mathsf{L}}
\newcommand{\sM}{\mathsf{M}}
\newcommand{\sN}{\mathsf{N}}
\newcommand{\sI}{\mathsf{I}}
\newcommand{\sS}{\mathsf{S}}
\newcommand{\sA}{\mathsf{A}}
\newcommand{\sQ}{\mathsf{Q}}

\newcommand{\ba}{\mathbf{i}}
\newcommand{\bb}{\mathbf{j}}
\newcommand{\bi}{\mathbf{i}}
\newcommand{\bj}{\mathbf{j}}

\newcommand{\conggr}{\cong_{\mathrm{gr}}}

\newcommand{\va}{\varphi}
\newcommand{\Ga}{\Gamma}
\newcommand{\ga}{\gamma}
\newcommand{\al}{\alpha}

\newcommand{\de}{\delta}
\newcommand{\la}{\lambda}

\newcommand{\Sg}{\Sigma_{\tau}}

\def\hsp{\mspace{1mu}}

\newcommand\SH{\operatorname{SH^0}}

\newcommand\SK[1][1]{\operatorname{SK}_{#1}}
\newcommand{\Nrd}[1][{}]{{\operatorname{Nrd}_{#1}}}

\newcommand{\GL}{\operatorname{GL}}
\newcommand{\EU}{\operatorname{EU}}
\newcommand{\Spin}{\operatorname{Spin}}
\newcommand{\SL}{\operatorname{SL}}

\newcommand{\Aut}{\operatorname{Aut}}

\newcommand{\chr}{\operatorname{char}}
\newcommand{\ind}{\operatorname{ind}}

\newcommand{\Gal}{\operatorname{Gal}}

\newcommand{\Br}{\operatorname{Br}}

\newcommand{\gr}{\operatorname{{\sf gr}}}
\newcommand{\im}{\operatorname{im}}

\theoremstyle{plain}
\newtheorem {lemma}{Lemma}[section]
\newtheorem {theorem}[lemma]{Theorem}
\newtheorem {thm}[lemma]{Theorem}
\newtheorem {corollary}[lemma]{Corollary}

\newtheorem {proposition}[lemma]{Proposition}
\newtheorem {prop}[lemma]{Proposition}

\theoremstyle{remark}
\newtheorem{remark}[lemma]{Remark}

\newtheorem {remarks}[lemma]{Remarks}
\newtheorem {example}[lemma]{Example}

\theoremstyle{definition}

\newtheorem{deff}[lemma]{Definition}

\numberwithin{equation}{section}

\begin{document}

\title[Unitary  $\SK$ of division algebras]
{Unitary  $\boldsymbol{\SK}$ of graded and valued division algebras, I}

\begin{thanks}
{The first author acknowledges the support of EPSRC first grant
scheme EP/D03695X/1. The second author would like to thank
the first author and Queen's University, Belfast for their
hospitality while the research for this paper was carried out.
}
\end{thanks}

\author{R. Hazrat}\address{
Department of Pure Mathematics\\
Queen's University\\
Belfast BT7 1NN\\
United Kingdom} \email{r.hazrat@qub.ac.uk}

\author{A. R. Wadsworth}
\address{
Department of Mathematics\\
University of California at San Diego\\
La Jolla, California 92093-0112\\
U.S.A.}
 \email{arwadsworth@ucsd.edu}

\begin{abstract}
The reduced unitary Whitehead group $\SK$ of a graded division
algebra equipped with a unitary involution (i.e.,  an involution of the second
kind) and graded by a torsion-free abelian group is studied.
It is shown that calculations in  the graded setting are much simpler than
their nongraded counterparts.   The bridge to the non-graded case is
 established by proving  that the unitary $\SK$ of a tame valued
division algebra wih a unitary involution over a henselian field
 coincides
with the unitary $\SK$ of its associated graded division
algebra.
As a consequence, the graded approach allows us not only to recover
results available in the literature with  substantially easier proofs,
but also to calculate the unitary $\SK$ for
much wider classes of  division algebras over henselian fields.

\end{abstract}

\maketitle

\section{Introduction}

Motivated by Platonov's striking work on  the reduced Whitehead group $\SK(D)$
of
valued division algebras $D$, see \cite{platonov,platsurvey},
V. Yanchevski\u\i, considered the unitary  analogue, $\SK(D, \tau)$, for a
division algebra~$D$ with unitary (i.e., second kind) involution $\tau$, see
\cite{yin,y,yinverse,yy}.
Working with division algebras over a field with henselian discrete (rank $1$)
valuation whose residue field also contains
a henselian discrete valuation, and carrying out formidable technical
calculations, he produced remarkable analogues to Platonov's results. By relating
$\SK(D,\tau)$ to data over the residue algebra, he showed not only that
  $\SK(D,\tau)$ could be nontrivial but that it could be any finite abelian group,
and he gave a formula in the bicyclic case expressing  $\SK(D,\tau)$
as a quotient of relative Brauer groups.
 Over the years since then   several approaches have been given  to
 understanding and calculating the (nonunitary) group $\SK$ using different methods,
notably by
Ershov~\cite{ershov}, Suslin~\cite{sus1,sus2}, Merkurjev and Rost~\cite{merk}
(For surveys on the group $\SK$, see \cite{platsurvey}, \cite{gille},
\cite{merk} or
\cite[\S6]{wadval}.) However, even after the passage of some  30 years,
there does not seem to have
been any improvement in calculating $\SK$ in  the unitary setting.
This may be due in part to the complexity of the formulas in Yanchevski\u\i's
work, and the difficulty in following some of his arguments.

This paper is a sequel to \cite{hazwadsworth} where the reduced Whitehead
group $\SK$ for a graded division algebra was studied. Here we consider the
reduced unitary Whitehead group of a graded division algebra
with unitary graded involution. As in our
previous work, we will see that the graded calculus is much easier and more
transparent than the non-graded one. We calculate  the unitary
$\SK$ in several important cases. We also show how this enables one to
calculate  the unitary $\SK$ of a tame
division algebra over a henselian field, by passage to the associated
graded division algebra.
The graded approach allows us not only to recover most of
Yanchevski\u\i's results in  \cite{y,yinverse, yy}, with very substantially
simplified proofs, but also extend them
to  arbitrary value groups and  to calculate the unitary $\SK$ for
wider classes of  division algebras.  There is a significant simplification
gained by considering arbitrary value groups from the outset, rather than
towers of discrete valuations. But the greatest gain comes from passage
to the graded setting, where the reduction  to arithmetic considerations
in  the degree $0$ division subring is quicker and more transparent.

We briefly describe our principal results. Let $\sE$ be a graded division
algebra, with torsion free abelian grade group $\Gamma_\sE$,
and let $\tau$ be a unitary
graded involution on $\sE$.
\lq\lq Unitary" means that
the action of $\tau$ on the center $\sT = Z(\sE)$ is nontrivial
(see~\S\ref{unitsk1}).
  The {\it reduced unitary Whitehead group}  for $\tau$ on~$\sE$ is defined
as
$$
\SK(\sE,\tau) \ = \ \big\{a\in \sE^*\mid \Nrd_\sE(a^{1-\tau})=1\big \}\big/
\big \langle a\in \sE^* \mid a^{1-\tau}=1 \big \rangle,
$$
where $\Nrd_\sE$ is the
reduced norm map $\Nrd_\sE\colon\sE^* \rightarrow \sT^*$
(see~\cite[\S3]{hazwadsworth}). Here, $a^{1-\tau}$ means $a\hsp\tau(a)^{-1}$.
Let $\sR=\sT^\tau = \{t\in \sT\mid \tau(t) = t\}\subsetneqq\sT$
(see~\S\ref{unitsk1}).
Let $\sE_0$ be the subring of homogeneous elements of degree~$0$ in $\sE$;
likewise for $\sT_0$ and $\sR_0$.
For an involution $\rho$ on $\sE_0$, $S _\rho(\sE_0)$ denotes
$\{a\in \sE_0 \mid \rho(a) = a\}$ and $\Sigma_\rho(\sE_0) = \langle
S_\rho(\sE_0)\cap \sE_0^*\rangle$.
Let $n$ be the index of $\sE$, and $e$  the exponent of the group
$\Gamma_\sE/\Gamma_\sT$. Since $[\sT:\sR]=2$, there are just two possible cases:
either  \ (i) $\sT$ is unramified over $\sR$, i.e.,
$\Gamma_\sT=\Gamma_\sR$; or \  (ii)~$\sT$~is~totally~ramified over~$\sR$,
i.e., $|\Gamma_\sT:\Gamma_\sR|=2$ . We will prove the following
formulas  for the unitary~$\SK$:

\begin{itemize}

\item[(i)] Suppose $\sT/\sR$ is unramified:

\medskip

\begin{itemize}
\item [$\bullet$] If $\sE/\sT$ is
unramified,  then $\SK(\sE,\tau)  \cong  \SK(\sE_0, \tau|_{\sE_0})$
(Prop.~\ref{unramified}).

\medskip

\item [$\bullet$] If $\sE/\sT$ is totally ramified,
then (Th.~\ref{sktotal}):
\begin{align*}
\SK(\sE,\tau)\ & \cong  \  \big \{a\in \sT_0^*\mid a^n\in \sR_0^*\}\big /
\{a\in \sT_0^*\mid a^e\in \sR_0^* \}
\\
&\cong  \  \big \{\omega \in \mu_n(\sT_0) \mid \tau(\omega)
=\omega^{-1}\big \}\big /\mu_e.
\end{align*}


\item [$\bullet$] If $\Gamma_\sE/\Gamma_\sT$ is cyclic,
and $\sigma$ is a generator of $\Gal(Z(E_0)/T_0)$,  then (Prop.~\ref{cyclic}):
\smallskip
\begin{itemize}
\item [$\circ$] $\SK(\sE,\tau) \ \cong  \  \{ a\in E_0^* \mid N_{Z(\sE_0)/\sT_0}
(\Nrd_{\sE_0}(a)) \in \sR_0  \}\big /\big(
\Sigma_\tau(\sE_0)\cdot \Sigma_{\sigma\tau}(\sE_0)\big)$.

\smallskip

 \item [$\circ$] If $\sE_0$ is a field, then $\SK(\sE,\tau)=1$.
\end{itemize}

 \medskip

\item[$\bullet$]  If $\sE$ has a maximal graded subfield
$\sM$ unramified over $\sT$ and another maximal graded subfield~$\sL$ totally
ramified over $\sT$, with $\tau(\sL ) =\sL$, then $\sE$ is semiramified and
(Cor.~\ref{seses})
\begin{equation*}
\SK(\sE,\tau) \ \cong \ \big\{a \in \sE_0 \mid N_{\sE_0/\sT_0}(a)\in \sR_0\big\}
 \, {\big/} \,  \textstyle{\prod\limits_{h\in \Gal(\sE_0/\sT_0)}}
\sE_0^{*h \tau}.
\end{equation*}
\end{itemize}


\item[(ii)] If $\sT/\sR$ is totally ramified,  then
$\SK(\sE,\tau)=1$ (Prop.~\ref{completely}).
\end{itemize}

The bridge between the graded and the non-graded henselian setting is
established by
Th.~\ref{involthm2}, which shows that for a tame  division
algebra $D$ over a henselian valued field with a unitary involution $\tau$,
${\SK(D,\tau)\cong \SK(\gr(D),\widetilde\tau)}$ where $\gr(D)$ is the
graded division algebra associated to $D$
by the valuation,
  and $\widetilde\tau$
is the graded involution on $\gr(D)$ induced by $\tau$ (see~\S\ref{unitary}).
Thus, each of the results listed above for graded division  algebras
yields analogous formulas for valued division algebras over a henselian
field, as illustrated in Example~\ref{toex} and Th.~\ref{appl}.
This recovers existing formulas, which were primarily for the
case with value group $\zz$ or $\zz\times \zz$, but with easier and
more transparent proofs than  those in the existing literature.
Additionally, our results apply   for any value groups whatever.
The especially simple case where $\sE/\sT$~is totally ramified
 and $\sT/\sR$~is unramified is entirely new.

In the sequel to this paper \cite{II}, the very interesting special
case will be treated where $\sE/\sT$ is semiramified (and $\sT/\sR$
is unramified) and $\Gal(\sE_0/\sT_0)$ is bicyclic.  This case
was the setting  of essentially all of Platonov's specifically computed
examples with nontrivial $\SK(D)$ \cite{platonov,plat76}, and likewise
Yanchevski\u\i's unitary examples in \cite{yinverse}
where the nontrivial $\SK(D, \tau)$ was fully computed.  This case is
not pursued here because it requires some more specialized arguments.
For such an $\sE$, it is known that $[\sE]$ decomposes
(nonuniquely) as  $[\sI\otimes _\sT \sN]$  in the
graded Brauer group of $\sT$, where $\sI$ is inertial over $\sT$ and
$\sN$  is nicely semiramified, i.e., semiramified and containing a
maximal  graded subfield totally ramified over $\sT$.  Then  a formula will
be given for $\SK(E)$ as a factor group of the relative Brauer
group $\Br(\sE_0/\sT_0)$ modulo other relative Brauer groups and the
class of $\sI_0$. An exactly analogous formula will be proved for
$\SK(\sE, \tau)$  in the unitary setting.

\section{Preliminaries}\label{prel}

Throughout this  paper we will be concerned with involutory division algebras and
involutory graded division algebras. In the non-graded setting, we will denote a
division algebra by $D$ and its  center by $K$;  this~$D$ is  equipped with an
involution
$\tau$, and we set  $F=K^\tau = \{a \in K\mid \tau(a) = a\}$.
 In the graded setting, we will write $\sE$ for a graded division algebra with
center $\sT$, and $\sR=\sT^\tau$ where $\tau$ is a graded involution on~$\sE$.
(This is consistent  with the notation used in~\cite{hazwadsworth}.)
Depending on
the context, we will write $\tau(a)$~or~$a^\tau$ for the action of the
involution on an element, and $K^\tau$ for the set of elements of $K$
invariant under $\tau$. Our convention is that $a^{\sigma\tau}$ means
$\sigma(\tau(a))$.

 In this section, we  recall  the notion of graded division algebras
and collect the  facts we need about them  in~\S\ref{pregda}. We will then
introduce the unitary and graded reduced unitary Whitehead groups
in~\S\ref{grinvols} and~\S\ref{unitsk1}.

\subsection{Graded division algebras}\label{pregda}

In this subsection we establish  notation and
recall some fundamental facts about graded division
algebras indexed by a totally ordered abelian group. For an extensive
study of such graded division algebras
and their relations with valued division algebras, we refer the reader
to~\cite{hwcor}. For generalities  on graded rings see~\cite{vanoy}.

Let
$\sR = \bigoplus_{ \ga \in \Gamma} \sR_{\ga}$ be a
graded ring, i.e.,
  $\Gamma$ is an abelian group,  and $\sR$ is a
unital ring such that each $\sR_{\ga}$ is a
subgroup of $(\sR, +)$ and
$\sR_{\ga} \cdot \sR_{\de} \subseteq \sR_{\ga +\de}$
for all $\ga, \de \in \Ga$. Set
\begin{itemize}
\item[] $\Gamma_\sR  \ = \  \{\ga \in \Ga \mid \sR_{\ga} \neq 0 \}$,
 \  the grade set of $\sR$;
\vskip0.05truein
\item[] $\sR^h  \ = \ \bigcup_{\ga \in
\Ga_{\sR}} \sR_{\ga}$,  \ the set of homogeneous elements of $\sR$.
\end{itemize}
For a homogeneous element of $\sR$  of degree $\gamma$, i.e., an
$r \in \sR_\gamma\mi\{0\}$, we write $\deg(r) = \gamma$.
Recall  that $\sR_{0}$ is a subring of $\sR$ and that for each $\ga \in
\Ga_{\sR}$, the group $\sR_{\ga}$ is a left and right $\sR_{0}$-module.
A subring $\sS$ of
$\sR$ is a \emph{graded subring} if $\sS= \bigoplus_{ \ga \in
\Gamma_{\sR}} (\sS \cap \sR_{\ga})$.  For example, the
center of $\sR$, denoted $Z(\sR)$, is a graded subring of
$\sR$.
If $\sT = \bigoplus_{ \ga \in
\Gamma} \sT_\gamma$ is another graded ring,
a {\it graded ring homomorphism} is a ring homomorphism
$f\colon \sR\to \sT$ with $f(\sR_\gamma) \subseteq \sT_\gamma$
for all $\gamma \in \Gamma$.  If $f$ is also bijective,
it is called a graded ring isomorphism;  we then write
$\sR\conggr \sT$.

For a graded ring $\sR$, a {\it graded left $\sR$-module} $\sM$ is
a left  $\sR$-module with a grading ${\sM=\bigoplus_{\ga \in \Ga'}
\sM_{\ga}}$,
where the $\sM_{\ga}$ are all abelian groups and $\Ga'$ is a
abelian group containing $\Ga$, such that $\sR_{\ga} \cdot
\sM_{\delta} \subseteq \sM_{\ga + \delta}$ for all $\ga \in \Ga_\sR,
\delta \in \Ga'$.
Then, $\Gamma_\sM$ and $\sM^h$ are
defined analogously to $\Gamma_\sR$ and~$\sR^h$.  We say that $\sM$ is
a {\it graded free} $\sR$-module if it has a base as a free
$\sR$-module consisting of homogeneous elements.

A graded ring $\sE = \bigoplus_{ \ga \in \Gamma} \sE_{\ga}$ is
called a \emph{graded division ring} if $\Ga$ is a
torsion-free abelian group and
every non-zero homogeneous
element of $\sE$ has a multiplicative inverse in $\sE$.
Note that the grade set  $\Ga_{\sE}$ is actually a group.
Also, $\sE_{0}$ is a division ring,
and  $\sE_\gamma$ is a $1$-dimensional
left and right $\sE_0$ vector space for every $\gamma\in \Gamma_\sE$.
Set $\sE_\gamma^* = \sE_\gamma \setminus\{0\}$.
The requirement that $\Gamma$ be torsion-free is made
because we are interested in graded division rings arising
from valuations on division rings, and all the grade groups
appearing there are torsion-free.  Recall that every
torsion-free abelian group $\Gamma$ admits total orderings compatible
with the group structure.  (For example, $\Gamma$ embeds in
$\Gamma \otimes _{\mathbb Z}\mathbb Q$ which can be given
a lexicographic total ordering using any base of it as a
$\mathbb Q$-vector space.)  By using any total ordering on
$\Gamma_\sE$, it is easy to see that $\sE$ has no zero divisors
and that $\sE^*$, the multiplicative group of units of $\sE$,
coincides with
 $\sE^{h} \mi \{0\}$ (cf. \cite[p.~78]{hwcor}).
Furthermore, the degree map
\begin{equation}\label{degmap}
\deg\colon \sE^* \rightarrow \Gamma_\sE
\end{equation} is a group epimorphism with kernel $\sE_0^*$.


By an easy adaptation of the ungraded arguments, one can see
that every graded  module~$\sM$  over a graded division ring
$\sE$ is graded free, and every two
homogenous bases have the same cardinality.
We thus call $\sM$ a \emph{graded vector space} over $\sE$ and
write $\dim_\sE(\sM)$ for the rank of~$\sM$ as a graded free $\sE$-module.
Let $\sS \subseteq \sE$ be a graded subring which is also a graded
division ring.  Then we can view $\sE$ as a graded left $\sS$-vector
space, and  we write $[\sE:\sS]$ for $\dim_\sS(\sE)$.  It is easy to
check the \lq\lq Fundamental Equality,"
\begin{equation}\label{fundeq}
[\sE:\sS] \ = \ [\sE_0:\sS_0] \, |\Gamma_\sE:\Gamma_\sS|,
\end{equation}
where $[\sE_0:\sS_0]$ is the dimension of $\sE_0$ as a left vector space
over the division ring $\sS_0$ and $|\Gamma_\sE:\Gamma_\sS|$ denotes the
index in the group $\Gamma_\sE$ of its subgroup $\Gamma_\sS$.

A \emph{graded field} $\sT$ is a commutative graded division ring.
Such a $\sT$ is an integral domain (as $\Gamma_\sT$ is torsion free),
so it has a quotient field,
which we denote $q(\sT)$.  It is known, see \cite[Cor.~1.3]{hwalg},
that $\sT$~is integrally closed in $q(\sT)$.  An extensive theory of
graded algebraic field extensions of graded fields has been developed in
\cite{hwalg}.

If $\sE$ is a graded division ring, then its center $Z(\sE)$ is clearly
a graded field.  {\it The graded division rings considered in
this paper will always be assumed finite-dimensional over their
centers.}  The finite-dimensionality assures that $\sE$
has a quotient division ring $q(\sE)$ obtained by central localization,
i.e., ${q(\sE) = \sE \otimes_\sT q(\sT)}$,  where $\sT = Z(\sE)$. Clearly,
$Z(q(\sE)) = q(\sT)$ and $\ind(\sE) = \ind(q(\sE))$, where
the index of $\sE$ is defined by $\ind(\sE)^2 = [\sE:\sT]$
(see \cite[p.~89]{hwcor}).
If $\sS$ is a graded field which is a graded subring of $Z(\sE)$
and $[\sE:\sS] <\infty$,
then $\sE$ is said to be a {\it graded division algebra} over~$\sS$.
We recall a fundamental connection between $\Gamma_\sE$ and $Z(\sE_0)$:
The field $Z(\sE_0)$
is Galois over $\sT_0$,  and there is a well-defined group epimorphism
\begin{equation}\label{surj}
\Theta_\sE\colon\Gamma_\sE \rightarrow \Gal(Z(\sE_0)/\sT_0)
\text{\,\,\,\,\, given by \,\,\,\,\, } \deg(e)\mapsto (z\mapsto eze^{-1}),
\end{equation}
for any $e\in \sE^*$.
(See \cite[Prop.~2.3]{hwcor} for a proof).

Let $\sE = \bigoplus_{\al \in \Ga_\sE} \sE_{\al}$ be a graded division
algebra with a graded center $\sT$ (with, as always, $\Ga_\sE$ a
torsion-free abelian group).
After fixing  some total ordering on $\Ga_\sE$,  define a function
$$
\la\colon  \sE\mi\{0\} \ra \sE^{\ast} \quad\textrm{ by }
\quad\la(\textstyle\sum c_{\ga}) =
c_{\de}, \text{ where $\de$ is minimal among the $ \ga \in \Ga_\sE$ with
$c_{\ga} \neq0$}.
$$
 Note that $\la(a) = a$ for $a \in \sE^{\ast}$, and
\begin{equation}\label{valhomin}
 \la(ab)= \la(a) \la(b) \textrm{ for all } a,b \in \sE\mi\{0\}.
\end{equation}

Let $Q = q(\sE)= \sE \otimes_\sT q(\sT)$, which is a division ring
as $\sE$ has no zero divisors and is finite-dimensional over $\sT$.
We can extend $\la$ to a map defined on all of $Q^{\ast} = Q\setminus\{0\}$
 as follows: for $q
\in Q^{\ast}$, write $q = ac^{-1}$  with $a \in \sE\mi\{0\}$, $c \in
Z(\sE)\mi\{0\}$, and set $\la(q) = \la(a) \la(c)^{-1}$. It follows
 from (\ref{valhomin}) that $\la\colon Q^{\ast} \ra \sE^{\ast}$ is
well-defined and is
a group homomorphism. Since the composition
\begin{equation}\label{inji}
\sE^{\ast} \hookrightarrow Q^{\ast}
\stackrel{\lambda}{\longrightarrow} \sE^{\ast}
\end{equation}
 is the identity, $\la$ is a
splitting map for the injection $\sE^{\ast} \hookrightarrow
Q^{\ast}$.

For a graded division algebra $\sE$ over its center $\sT$, there is
 a reduced norm map $\Nrd_\sE\colon\sE^*\rightarrow \sT^*$
(see~\cite[\S3]{hazwadsworth}) such that  for $a\in \sE$ one has
$\Nrd_\sE(a)=\Nrd_{q(\sE)}(a)$. The {\it reduced Whitehead group},
$\SK(\sE)$, is defined as
$\sE^{(1)}/\sE'$, where $\sE^{(1)}$ denotes
the set of elements of $\sE^*$ with  reduced norm 1, and $\sE'$ is the
commutator subgroup $[\sE^*,\sE^*]$ of $\sE^*$. This group was studied in detail
in~\cite{hazwadsworth}.
We will be using the following facts which were established in that paper:
\begin{remarks}\label{grfacts} Let $n = \ind(\sE)$.
\begin{enumerate}[\upshape(i)]
\item For $\gamma\in \Gamma_\sE$, if
$a\in \sE_\gamma$ then $\Nrd_\sE(a)\in \sE_{n\gamma}$.  In particular,
$\sE^{(1)}$ is a subset of $\sE_0$.

\smallskip

\item  If $\sS$ is any graded subfield of $\sE$
containing $\sT$ and $a\in \sS$, then
$\Nrd_\sE(a) =  N_{\sS/\sT}(a)^{n/[\sS:\sT]}$.

\smallskip

\item \label{rnrd} Set
\begin{equation}\label{dlambda}
\dlambda  \ = \  \ind(\sE)\big/\big(\ind(\sE_0) \,
[Z(\sE_0) \!: \!\sT_0]\big).
\end{equation}
  If $a\in \sE_0$, then,
\begin{equation}\label{nrddo}
\Nrd_\sE(a)  \ = \   N_{Z(\sE_0)/\sT_0}\Nrd_{\sE_0}(a)^
 {\,\dlambda}
\in \sT_0.
\end{equation}

\smallskip

\item \label{normal} If $N$ is a normal subgroup of
$\sE^{\ast}$, then $N^{n} \subseteq
\Nrd_\sE(N)[\sE^{\ast}, N]$.

For  proofs of (i)-(iv) see \cite[Prop.~3.2 and~3.3]{hazwadsworth}.

\smallskip

\item \label{torik} $\SK(\sE)$ is $n$-torsion.
\begin{proof}
By taking  $N=\sE^{(1)}$, the assertion follows from (\ref{normal}).
\end{proof}
\end{enumerate}
\end{remarks}

A graded division algebra $\sE$ with center $\sT$ is
said to be {\it inertial} (or {\it unramified}) if $\Gamma_\sE=\Gamma_\sT$.
From~(\ref{fundeq}), it then follows  that $[\sE:\sT]=[\sE_0:\sT_0]$;
indeed, $\sE_0$ is central simple over $\sT_0$ and $\sE \cong_{\gr} \sE_0 \otimes _{\sT_0} \sT$.
At the other extreme, $\sE$ is said to be {\it totally ramified}
if $\sE_0=\sT_0$. In an intermediate case  $\sE$ is
said to be {\it semiramified} if
 $\sE_0$ is a field and ${[\sE_0:\sT_0]=|\Gamma_\sE:\Gamma_\sT|=\ind(\sE)}$.
These definitions are motivated by  analogous  definitions for
valued division algebras (\cite{wadval}).
Indeed,  if a tame valued division algebra
is unramified, semiramified, or totally ramified, then so is its
associated graded division algebra.
Likewise, a  graded field extension $\sL$~of~
$\sT$ is said to be  {\it inertial} (or {\it unramified})
if $\sL\cong_{\gr} \sL_0 \otimes_{\sT_0}\sT$ and the field $\sL_0$
is separable over $\sT_0$. At the other extreme,
$\sL$ is {\it totally ramified} over $\sT$ if
$[\sL:\sT] = |\Gamma_\sL:\Gamma_\sT|$.  A graded division algebra
$\sE$ is said to be {\it inertially split} if $\sE$ has a maximal
graded subfield $\sL$ with $\sL$ inertial  over $\sT$.
When this occurs, $\sE_0 = \sL_0$,  and $\ind(\sE)= \ind(\sE_0) \,
[Z(\sE_0) \!: \!\sT_0]$ by Lemma~\ref{dlambdafacts} below.
In particular, if $\sE$ is  semiramified then $\sE$ is inertially
split,
$\sE_0$ is
abelian Galois over $\sT_0$, and the canonical map
$\Theta_\sE\colon \Gamma_\sE\rightarrow \Gal(\sE_0/\sT_0)$ has kernel
$\Gamma_\sT$  (see~\eqref{surj} above and
~\cite[Prop.~2.3]{hwcor}).

\begin{lemma}\label{dlambdafacts}
Let  $\sE$ be a graded division algebra  with center $\sT$.
 For the  $\dlambda$ of~\eqref{dlambda},
$\dlambda^2 = |\ker(\Theta_\sE)/\Gamma_\sT|$.  Also,
$\dlambda = 1$ iff $\sE$ is inertially split.
\end{lemma}

\begin{proof}
Since $\Theta_\sE$ is surjective,  $\Gamma_\sT\subseteq\ker(\Theta_\sE)$,
and $Z(\sE_0)$ is Galois over $\sT_0$, we have
\begin{align*}
\dlambda^2  \ &= \ \ind(\sE)^2 \, \big/ \,
\big(\ind(\sE_0)^2 \, [Z(\sE_0):\sT_0]^2\big)
 \ = \ [\sE:\sT] \, \big/\big([\sE_0:Z(\sT_0)] \, [Z(\sE_0):\sT_0] \,
|\Gal(Z(\sE_0)/\sT_0)|\big)\\
& =  \ [\sE_0:\sT_0] \, |\Gamma_\sE/\Gamma_\sT| \, \big/ \,
\big([\sE_0:\sT_0] \, |\im(\Theta_\sE)|\big) \ = \ |\ker(\Theta_\sE)/\Gamma_\sT|.
\end{align*}

Now, suppose $M$ is a maximal subfield of $\sE_0$
with $M$ separable over $\sT_0$.  Then,
$M\supseteq Z(\sE_0)$ and $[M:Z(\sE_0)] = \ind(\sE_0)$.  Let
$\sL = M\otimes _{\sT_0} \sT$ which is a graded subfield of $\sE$
inertial over $\sT$, with $\sL_0 = M$.  Then,
$$
[\sL:\sT] \ = \ [\sL_0:\sT_0] \ = \ [\sL_0:Z(\sE_0)]\, [Z(\sE_0):\sT_0]
 \ = \ \ind(\sE)/\dlambda.
$$
Thus, if $\dlambda = 1$, then $\sE$ is inertially
split, since $\sL$ is a maximal graded subfield of $\sE$ which is
inertial over $\sT$.  Conversely, suppose $\sE$ is inertially
split, say $\sI$ is a maximal graded subfield of $\sE$ with
$\sI$  inertial over $\sT$.  So, $[\sI_0:\sT_0] = [\sI:\sT] =
\ind(\sE)$. Since $\sI_0Z(\sE_0)$ is a subfield of $\sE_0$,
we have
\begin{align*}
[\sI_0:\sT_0] \ &\le \ [\sI_0Z(\sE_0) :\sT_0] \ = \
[\sI_0 Z(\sE_0):Z(\sE_0)]\, [Z(\sE_0):\sT_0] \\
&\le  \
\ind (\sE_0)\, [Z(\sE_0):\sT_0] \ = \ \ind(\sE)/\dlambda
 \ = \ [\sI_0:\sT_0] /\dlambda.
\end{align*}
So, as $\dlambda$ is a positive integer, $\dlambda = 1$.
\end{proof}

\subsection{Unitary $\SK$ of division algebras}\label{grinvols}

We begin with a description of unitary $K_1$ and $\SK$ for a division
algebra with an involution. The analogous definitions for graded
division algebras will be given in ~\S\ref{unitsk1}.

Let $D$ be a division ring finite-dimensional over its center $K$ of
index $n$, and let $\tau$ be an involution on~$D$,  i.e., $\tau$ is
an antiautomorphism of $D$ with $\tau^2 =\id$.
Let
\begin{itemize}
\item[] $\quad S_\tau(D)  \, =  \, \{d \in D\mid\tau(d) = d \};$
\vskip0.05truein
\item [] $\quad \Sigma_{\tau}(D)  \, = \,
\left< S_{\tau}(D) \cap D^{\ast} \right>.$
\end{itemize}

Note that $\Sigma_{\tau}(D)$ is a normal subgroup of $D^{\ast}$. For,
if $a \in S_{\tau}(D)$, $a \neq 0$, and $b \in D^{\ast}$, then
$bab^{-1} = [ba\tau(b)][b\tau(b)]^{-1}\in \Sigma_\tau(D)$,
as $ba\tau(b), b\tau(b) \in S_\tau(D)$.

Let $\va$ be an isotropic $m$-dimensional, nondegenerate skew-Hermitian form
over $D$
with respect to an involution $\tau$ on $D$.  Let
$\rho$ be the involution on $M_m(D)$ adjoint to $\va$, let
$U_m(D) = \{a \in M_m(D)\mid a\rho(a) = 1\}$ be the unitary group
associated to $\va$,  and let $\EU_m(D)$ denote
the normal subgroup of $U_m(D)$ generated by the unitary transvections.
For  $m>2$,
the Wall spinor norm map
$\Theta\colon U_m(D) \rightarrow D^*/\Sigma_\tau(D)D'$ was developed
in \cite{wall},  where it was shown that  $\ker(\Theta) = \EU_m(D)$.
Here, $D'$ denotes the multiplicative commutator group $[D^*,D^*]$.
Combining this
with~\cite[Cor.~1 of~\S 22]{draxl} one obtains  the commutative diagram:
\begin{equation}\label{unnrd}
\begin{split}
\xymatrix{U_m(D)\big/\EU_m(D) \ar[r]_-{\cong}^-\Theta  \ar[d] &
D^*\big/\big ( \Sigma_\tau(D)D' \big )\ar[d]^{1-\tau}\\
\GL_m(D)\big/E_m(D) \ar[r]^-{\det} \ar[d]_{\Nrd} & D^*\big/D'
\ar[d]^{\Nrd}\\
K^* \ar[r]^-{\id} & K^*}
\end{split}
\end{equation}
where the map $\det$ is the Dieudonn\'e determinant and
$1-\tau\colon D^*/\big ( \Sigma_\tau(D)D' \big ) \longrightarrow D^*/D'$
is defined as
$x\Sigma_\tau(D)D'\mapsto x^{1-\tau}D'$, where $x^{1-\tau}$ means $x\hsp\tau(x)^{-1}$
(see~\cite[6.4.3]{hahn}).

From the diagram, and parallel to the ``absolute'' case, one defines
the {\it unitary Whitehead group},
$$
K_1(D,\tau)\ = \ D^*/\big (\Sigma_\tau(D)D'\big ).
$$

For any  involution $\tau$ on $D$, recall that
\begin{equation}\label{taunrd}
\Nrd_D(\tau(d))  \ = \  \tau(\Nrd_D(d)),
\end{equation} for any  $d\in D$. For, if $p \in K[x]$ is the
minimal polynomial of $d$ over $K$,
then $\tau(p)$ is the minimal polynomial of $\tau(d)$ over $K$
(see also~\cite[\S22, Lemma~5]{draxl}).

 We consider two cases:

\subsubsection{ Involutions of the first kind}\label{firstk}
In this case the center $K$ of $D$ is elementwise invariant under the
involution, i.e.,
$K \subseteq S_\tau(D)$. Then $S_\tau(D)$ is a $K$-vector space.  The
 involutions of this kind are further subdivided into two types:
{\it orthogonal} and
{\it symplectic} involutions ~(see \cite[Def.~2.5]{kmrt}).
By ~(\cite[Prop.~2.6]{kmrt}), if $\chr(K)\not = 2$ and $\tau$ is orthogonal
then, $\dim_K(S_\tau(D))=n(n+1)/2$, while if $\tau$ is symplectic then
$\dim_K(S_\tau(D))=n(n-1)/2$. However, if $\chr(K)=2$, then
$\dim_K(S_\tau(D))=n(n+1)/2$ for each type.

If $\dim_K(S_\tau(D))=n(n+1)/2$, then  for any
$x \in D^*$, we have $xS_\tau(D) \cap S_\tau(D) \not = 0$
by dimension count; it then follows that
$D^*=\Sigma_\tau(D)$, and thus $K_1(D,\tau)=1$. However,
in the case $\dim_K(S_\tau(D))=n(n-1)/2$, Platonov showed that $K_1(D,\tau)$ is not in general
trivial, settling  Dieudonn\'e's conjecture in negative~\cite{pld}.
Note that whenever $\tau$  is of the first kind we have
$\Nrd_D(\tau(d)) = \Nrd_D(d)$ for all $d\in D$, by \eqref{taunrd}.
Thus,  $K_1(D,\tau)$ is sent to  the identity under
the composition   $\Nrd \circ (1-\tau)$. This explains why
 one does not consider the kernel of this map, i.e., the
unitary $\SK$, for involutions of the first kind.
If $\chr(K) \not =2$ and
$\tau$~is symplectic,
then as the $m$-dimensional form $\va$ over $D$ is skew-Hermitian,
its associated  adjoint involution~$\rho$  on $M_m(D)$ is
of orthogonal type, so there is  an associated  spin group
$\Spin(M_m(D), \rho)$.
For any $a\in S(D)$ one then has $\Nrd_D(a) \in K^{*2}$
(\cite[Lemma~2.9]{kmrt}).
One defines  $K_1\Spin(D,\tau)=R(D)/\big(\Sigma_\tau(D)D'\big)$,
where
$R(D)=\{d\in D^* \mid \Nrd_D(d) \in K^{*2}\}$. This group
is  related to $\Spin(M_m(D), \rho)$, and has been studied
in~\cite{monyan},  parallel to the work on absolute  $\SK$ groups and
unitary $\SK$ groups for unitary involutions.

 \subsubsection{Involutions of the second kind $($unitary involutions$)$} \label{secondk}
 In this case $K \not \subseteq S_\tau(D)$.
 Then, let \break ${F=K^{\tau}\ (=K\cap
S_\tau(D)\,)}$, which is a subfield of $K$ with $[K:F] = 2$.
 It was already observed by Dieudonn\'e that
$U_m(D) \not = \EU_m(D)$.  An important property proved
 by Platonov and Yanchevski\u\i, which we will use frequently,
is that
\begin{equation}\label{primeinsigma}
D' \,  \subseteq  \, \Sigma_\tau(D).
\end{equation}
(For a proof, see~\cite[Prop.~17.26]{kmrt}.) Thus
$K_1(D,\tau)=D^*/\Sigma_\tau(D)$, which  is not trivial in general.
The kernel of the
map $\Nrd \circ (1-\tau)$ in  diagram~\eqref{unnrd}, is called the
{\it reduced unitary Whitehead group}, and denoted by $\SK(D,\tau)$.
Using~(\ref{taunrd}), it is straightforward to see that
$$
\SK(D,\tau) \ = \  \Sigma'_{\tau}(D) / \Sigma_{\tau}(D), \quad \text{where}
\quad\Sigma'_{\tau}(D)  \ = \  \{ a \in D^{\ast} \mid \Nrd_D(a) \in F^{\ast}\}.
$$
Note that  we use the notation $\SK(D,\tau)$ for the reduced unitary
 Whitehead group as opposed to Draxl's notation
$\operatorname{USK_1}(D,\tau)$ in~\cite[p.~172]{draxl} and
Yanchevski\u\i's notation
$\operatorname{SUK_1}(D,\tau)$~\cite{y} and the notation
$\operatorname{USK_1}(D)$ in \cite{kmrt}.

Before we define the corresponding groups in the graded setting, let us
recall that all the groups above fit in Tits' framework
\cite{tits}
of the
{\it Whitehead group} $W(G,K)=G_K/G_K^+$ where $G$ is an almost simple,
simply connected linear algebraic group defined over an infinite field
$K$, with $\chr(K)\not =2$, and $G$ is isotropic over $K$. Here, $G_K$ is the
set of $K$-rational points of $G$ and $G_K^+$, is the subgroup of $G_K$,
generated by the unipotent radicals of the minimal $K$-parabolic
subgroups of $G$. In this setting, for $G_K=\SL_n(D)$, $n>1$, we have
$W(G,K)=\SK(D)$; for $\tau$ an involution of first or second  kind
on $D$ and $F = K^\tau$, for  $G_F=\SL_n(D,\tau):=\SL_n(D)\cap U_n(D)$
we have  $W(G,F)=\SK(D,\tau)$; and for $\tau$ a symplectic involution
on~ $D$ and $\rho$ the adjoint involution of an $m$-dimensional
isotropic skew-Hermitian
 form over $D$ with $m\ge 3$,
for the spinor group
$G_K=\Spin(M_m(D),\rho)$ we have $W(G,K)$
is a double cover of $K_1 \Spin(D,\tau)$ (see~\cite{monyan}).

\subsection{Unitary $\SK$ of graded division algebras}\label{unitsk1}
We will now introduce the unitary $K_1$ and $\SK$ in the graded setting.
Let $\sE = \bigoplus_{\ga \in \Ga_\sE} \sE_{\ga}$ be a graded division
ring (with $\Ga_\sE$ a torsion-free abelian group) such that
$\sE$ has finite dimension $n^2$ over its center $\sT$, a graded field.
Let $\tau$ be a graded involution of~$\sE$, i.e., $\tau$~is
an antiautomorphism of $\sE$ with $\tau^2 =\id$ and $\tau(\sE_{\ga}) =
\sE_{\ga}$ for each $\ga \in \Ga_\sE$. We define $S_{\tau}(\sE)$ and
$\Sigma_{\tau}(\sE)$, analogously to the non-graded cases, as the set of
elements of $\sE$ which are invariant under~$\tau$, and the
multiplicative group
generated by the nonzero homogenous elements of $S_\tau(\sE)$, respectively.
We say  the involution of the {\it first kind} if all the elements of the
center $\sT$ are invariant under $\tau$; it is of the {\it second kind}
(or {\it unitary})
otherwise. If $\tau$ is of the first kind then,  parallel to the
non-graded case, either
$\dim_\sT(S_\tau(\sE))=n(n+1)/2$ or
$\dim_\sT(S_\tau(\sE))=n(n-1)/2$.
Indeed,  one can show these equalities
by arguments analogous to the
nongraded case as in the proof of \cite [Prop.~2.6(1)]{kmrt}, as $\sE$ is split 
by a graded maximal subfield and the Skolem--Noether theorem is available
in the graded setting (\cite[Prop.~1.6]{hwcor}).  (These equalities  can also
be obtained
by passing to the quotient division algebra as is done in
Lemma~\ref{invfacts}\eqref{type} below.)


 Define the {\it
 unitary Whitehead group}
$$
K_1(\sE,\tau) \ = \ \sE^*/\big(\Sigma_\tau(\sE)\sE'\big),
$$
where $\sE' = [\sE^*,\sE^*]$.
If $\tau$ is of the first kind,
$\chr(\sT)\not =2$, and $\dim_\sT(S_\tau(\sE))=n(n-1)/2$, a proof similar
to \cite[Prop.~2.9]{kmrt}, shows that if $a
\in S_\tau(\sE)$ is homogeneous, then $\Nrd_\sE(a) \in \sT^{*2}$
(This can also be verified by passing to the quotient division
algebra, then using Lemma~\ref{invfacts}\eqref{type} below and invoking the
corresponding result for ungraded division algebras.)   For
this type of involution, define the {\it spinor Whitehead group}
$$
K_1\Spin(\sE,\tau) \ = \ \{a \in \sE^*\mid
\Nrd_\sE(a) \in \sT^{*^2} \} \, / \big( \Sigma_\tau(\sE)\sE'\big).
$$

When the graded involution $\tau$ on $\sE$ is unitary, i.e.,  $\tau|_\sT \neq
\id$, let $\sR= \sT^{\tau}$, which is a graded subfield of~$\sT$
with $[\sT:\sR]=2$. Furthermore, $\sT$ is Galois over $\sR$, with
$\Gal (\sT/\sR) = \{\id, \tau|_\sT\}$. (See \cite{hwalg}
for Galois theory for graded field extensions.) Define the
{\it reduced unitary Whitehead group}
\begin{equation}\label{skdef}
\SK (\sE, \tau)  \ = \  \Sigma'_{\tau}(\sE) \, / \,
(\Sigma_{\tau}(\sE) \, \sE')
 \ = \ \Sigma'_{\tau}(\sE) \, / \,
\Sigma_{\tau}(\sE) ,
\end{equation}
where
$$
\Sigma'_{\tau}(\sE) \ = \ \big\{a\in \sE^*\mid
\Nrd_\sE(a^{1-\tau})=1\big \}  \ = \ \{a\in \sE^*\mid \Nrd_\sE(a) \in \sR^*\}
$$
and
$$
\Sigma_{\tau}(\sE) \ = \ \langle a\in \sE^* \mid
a^{1-\tau}=1 \big \rangle \ = \ \langle S_\tau(\sE) \cap \sE^*\rangle.
$$
Here, $a^{1-\tau}$ means $a\tau(a)^{-1}$. See Lemma~\ref{invfacts}(iv)
below for the second equality in \eqref{skdef}.
The group $\SK(\sE, \tau)$ will be the main focus of the rest of the paper.

We will use the following facts repeatedly:

\begin{lemma}\label{invfacts}\hfill
\begin{enumerate}[\upshape(i)]
\item \label{type}
Any graded involution on $\sE$ extends uniquely to
an involution of the same kind
$($and type$)$  on $Q=q(\sE)$.

\smallskip

\item \label{sigma} For any graded involution $\tau$ on $\sE$, and its
extension to $Q=q(\sE)$, we have $\Sigma_\tau(Q) \cap \sE^*
\subseteq \Sigma_\tau(\sE)$.

\smallskip

\item \label{firstsig} If $\tau$ is a graded involution of
the first kind on $\sE$ with $\dim_\sT(S_\tau(\sE))=n(n+1)/2$,
 then $\Sigma_\tau(\sE)=\sE^*$.

\smallskip

\item \label{yan} If $\tau$ is a unitary graded involution on $\sE$,
then $\sE'\subseteq \Sigma_\tau(\sE)$.

\smallskip

\item \label{tori} If $\tau$ is a unitary graded involution on $\sE$,
then $\SK(\sE,\tau)$ is a torsion group of bounded exponent
dividing  $n=\ind(\sE)$.

\end{enumerate}\end{lemma}

\begin{proof}\hfill

(i) Let $\tau$ be a graded involution on $\sE$.  Then
$q(\sE)=\sE\otimes_\sT q(\sT)=
\sE\otimes_\sT (\sT\otimes_{\sT^\tau}q(\sT^\tau))=
\sE\otimes_{\sT^\tau} q(\sT^\tau)$. The unique extension
of $\tau$  to $q(\sE)$ is $\tau  \otimes \id_{q(\sT^\tau)}$,
which we denote simply as $\tau$.
It then follows that
$S_\tau(q(\sE))=S_\tau(\sE)\otimes_{\sT^\tau} q(\sT^\tau)$.
Since $q(\sT^\tau)=q(\sT)^\tau$, the assertion follows.

(ii) Note that for the
 map $\lambda$ in the sequence (\ref{inji}) we  have
$\tau(\lambda(a))=\lambda(\tau(a))$ for all $a\in Q^*$.
Hence, $\lambda(\Sigma_\tau(Q)) \subseteq\Sigma_\tau(\sE)$.
Since $\lambda|_{\sE^*}$ is the identity,   we have
$\Sigma_\tau(Q) \cap \sE^* \subseteq \Sigma_\tau(\sE)$.

(iii) The extension of the graded  involution $\tau$  to
$Q=q(\sE)$, also denoted $\tau$, is of the first kind with
$\dim_Q(S_\tau(Q))=n(n+1)/2$ by (\ref{type}). Therefore
$\Sigma_\tau(Q)=Q^*$ (see~\S\ref{firstk}). Using
(\ref{sigma}) now, the assertion follows.

(iv) Since $\tau$ is a unitary graded involution,   its extension to
$Q=q(\sE)$ is  also  unitary,  by (\ref{type}). But
${Q' \subseteq \Sigma_\tau(Q)}$,
as noted in \eqref{primeinsigma}.  From~(\ref{inji}) it follows that
$Q' \cap \sE^* = \sE'$. Hence, using \eqref{sigma}, \break
${\sE' \subseteq \sE^* \cap Q' \subseteq \sE^* \cap \Sigma_\tau(Q)
\subseteq \Sigma_\tau(\sE)}$.

(v) Setting $N=\Sigma_\tau'(\sE)$,
Remark~\ref{grfacts}(\ref{normal}) above, coupled with the fact
that $\sE' \subseteq \Sigma_\tau(\sE)$ (\ref{yan}), implies that
$\SK(\sE,\tau)$ is an $n$-torsion group. This assertion also follows
by using (\ref{sigma}) which implies the natural map $\SK(\sE,\tau)
\to \SK(Q,\tau)$ is injective and the fact that unitary
$\SK$ of a division algebra of index $n$ is
$n$-torsion~(\cite[Cor. to~2.5]{y}).
\end{proof}

\subsection{Generalized dihedral groups and field extensions}

The nontrivial case of  $\SK(\sE,\tau)$ for $\tau$ a unitary
graded involution turns out to be when
 $\sT = Z(\sE)$ is unramified over $\sR = \sT^\tau$ (see \S\ref{unram}).
When that occurs, we will see in
Lemma~\ref{unramfacts}\eqref{eight} below that $Z(\sE_0)$ is a so-called
generalized dihedral extension over $\sR_0$.
We now give the definition and observe a few easy
 facts about generalized dihedral groups and extensions.

\begin{deff}\label{gendi}\hfill
\begin{enumerate}[\upshape(i)]
\item \label{deffo1} A group $G$ is said to be {\it generalized dihedral} if $G$ has a
subgroup $H$ such that $[G:H]=2$ and every $\tau \in G \backslash H$
satisfies $\tau^2=\id$.

Note that if $G$ is generalized dihedral and $H$ the distinguished subgroup,
then $H$ is abelian and $(h\tau)^2=\id$, for all $\tau \in G\backslash H$
and $h\in H$. Thus, $\tau^2=\id$ and $\tau h \tau ^{-1}=h^{-1}$ for all
$\tau \in G\backslash H, h \in H$. Furthermore, every subgroup of $H$ is
normal in $G$. Clearly every dihedral group is generalized dihedral, as
is every elementary abelian $2$-group. More generally, if $H$ is any
abelian group and $\chi \in \Aut(H)$ is the map $h \mapsto h^{-1}$, then
the semi-direct product
$H \rtimes_i \langle \chi \rangle$ is a generalized dihedral group,
where $i\colon
\langle \chi \rangle \rightarrow \Aut(H)$ is the inclusion map.
It is easy to check that every generalized dihedral group is isomorphic to
such a semi-direct product.

\item \label{deffo2} Let $F \subseteq K \subseteq L$ be fields with $[L:F]<\infty$ and
$[K:F]=2$. We say that $L$ is {\it generalized dihedral for $K/F$} if $L$ is
Galois over $F$ and every element of $\Gal(L/F) \backslash \Gal(L/K)$
has order $2$, i.e., $\Gal(L/F)$ is a generalized dihedral group.
Note that when this occurs, $L$ is compositum of fields $L_i$
 containing~$K$ with each
$L_i$ generalized dihedral for $K/F$ with $\Gal(L_i/K)$ cyclic, i.e.,
$L_i$ is Galois over $F$ with $\Gal(L_i/F)$ dihedral (or a Klein $4$-group).
Conversely, if $L$ and $M$ are generalized dihedral for $K/F$ then so is
their compositum.

\end{enumerate}
\end{deff}
\begin{example} Let $n\in \mathbb N$, $n \geq 3$, and let $F\subseteq K$ be fields with
$[K:F]=2$ and $K=F(\omega)$, where $\omega$ is a primitive $n$-th root of unity
(so $\chr(F) \nmid n$). Suppose the non-identity element of $\Gal(K/F)$ maps
$\omega$ to~$\omega^{-1}$. For any $c_1,\dots,c_k \in F^*$, if
$\omega \not \in F(\sqrt[n]{c_1},\dots,\sqrt[n]{c_k})$, then
$K(\sqrt[n]{c_1},\dots,\sqrt[n]{c_k})$ is generalized dihedral for $K/F$.
\end{example}

\section{Henselian to graded reduction}\label{unitary}

The main goal of this section is to prove an isomorphism
between  the unitary $\SK$ of a valued division
algebra with involution over a henselian field
  and the graded $\SK$ of  its associated graded division algebra.
We first recall how to associate a graded division algebra to a valued
division algebra.

Let $D$ be a division algebra finite dimensional over its
center $K$, with a valuation
$v\colon  D^{\ast} \ra
\Ga$. So, $\Ga$ is a totally ordered abelian group,
and $v$ satisifies the conditions that for all
$a, b \in D^{\ast}$,
\begin{enumerate}
\item $ \qquad \, v(ab)  \, = \,  v(a) + v(b)$;

\item $\quad v(a+b)  \, \geq \,  \min \{v(a),v(b) \}\;\;\;\;\; (b \neq -a).$
\end{enumerate}
Let
\begin{align*}
V_D  \ &= \  \{ a \in D^{\ast} \mid v(a) \geq 0 \}\cup\{0\},
\text{ the
valuation ring of $v$};\\
M_D  \ &= \  \{ a \in D^{\ast} \mid v(a) > 0
\}\cup\{0\}, \text{ the unique maximal left (and right) ideal
 of $V_D$}; \\
\overline{D}  \ &= \  V_D / M_D, \text{ the residue
division ring of $v$ on $D$; and} \\
\Ga_D  \ &= \  \mathrm{im}(v), \text{ the value
group of the valuation}.
\end{align*}

Now let $K$ be a field with a valuation $v$,
and suppose $v$ is {\it henselian}; that is,
  $v$ has a unique extension to every  algebraic
field extension of $K$. Recall that a
 field extension $L$ of $K$ of degree~$n<\infty$ is
said to be {\it tamely ramified} or {\it tame} over $K$
if, with respect to the unique extension of $v$ to $L$,
the residue field
$\overline L$ is separable over
$\overline K$ and $\chr({\overline K})\nmid
\big ( n\big/[\overline L:\overline K] \big )$. Such an $L$~is
necessarily {\it defectless} over $K$,
i.e., $[L:K] = [\overline L:\overline K] \, |\Gamma_L:\Gamma_K|$,
by \cite[Th.~3.3.3]{EP} (applied to $N/K$ and $N/L$, where
$N$ is a normal closure of $L$ over $K$).
Along the same lines, let $D$ be a
division algebra with center $K$ (so, by convention,
$[D:K] < \infty$); then the henselian valuation  $v$ on $K$ extends uniquely
to a valuation
on $D$ ~(\cite{wad87}). With respect to this valuation, $D$
is said to be
{\it tamely ramified} or {\it tame} if the center $Z(\overline D)$ is
 separable over $\overline K$ and ${\chr({\overline K}) \nmid
\big[\ind(D)\big/\big(\ind(\overline D)[Z(\overline D):\overline K]\big)}\big]$.
Recall from \cite[Prop.~1.7]{jw}, that whenever the field extension
$Z(\ov D)/\ov K$ is separable, it is abelian Galois.
It is known  that
$D$ is tame if and only if $D$~is split by the
maximal tamely ramified field extension of~$K$, if and only if
$\chr(\ov K) = 0$
or $\chr(\ov K) = p\ne 0$ and the $p$-primary component of~
$D$ is inertially split, i.e., split by the maximal unramified
extension of~$K$~(\cite[Lemma~6.1]{jw}).
We say $D$ is \emph{strongly tame}
if $\chr(\overline K)\nmid\ind(D)$.
Note that strong tameness implies tameness.
This is clear from the last characterization of tameness,
or from \eqref{Ostrowski} below.
 Recall also from \cite[Th.~3]{M}, that for a valued
division algebra $D$ finite dimensional over its  center $K$ (here not necessarily
henselian), we have the \lq\lq Ostrowski theorem"
\begin{equation}\label{Ostrowski}
[D:K]
\ = \ q^k\,[\overline D:\overline K] \,|\Gamma_D:\Gamma_K|,
\end{equation}
where $q=\chr({\overline D})$ and $k \in \mathbb Z$ with $k\geq 0$
(and  $q^k = 1$ if $\chr(\overline D) = 0$).
If $q^k = 1$ in equation~   \eqref{Ostrowski}, then $D$~is
said to be  {\it defectless} over $K$.
For background on  valued division algebras,
see~\cite{jw} or the survey paper~\cite{wadval}.

\begin{remark}\label{shensel}
If a field $K$ has a henselian valuation $v$ and $L$ is a subfield of
$K$ with $[K:L] <\infty$, then the restriction $w = v|_L$ need not be
henselian.  But it is easy to see that $w$ is then \lq\lq semihenselian,"
i.e., $w$~has more than one but only finitely many different extensions to
a separable closure $L_{\textrm{sep}}$ of $L$.
See \cite{engler} for a thorough analysis of
semihenselian valuations. Notably, Engler shows
that  $w$~is semihenselian iff the residue
field $\ov L_w$ is algebraically closed but
there is a henselian valuation $u$ on $L$ such that $u$
 is a proper coarsening of $w$ and  the residue  field $\ov L_u$ is
real closed. When this occurs, $\chr(L) = 0$, $L$ is formally real,
$w$ has exactly two extensions to
$L_{\textrm{sep}}$, the value group $\Gamma_{L,w}$
has a nontrivial divisible subgroup, and the henselization
of $L$ re $w$ is $L(\sqrt{-1})$, which lies in $K$.
For example, if we take any prime number $p$,
let $w_p$ be the
$p$-adic discrete valuation on $\qq$, and let $L = \{r\in \rr\,|\ r
\text{ is algebraic over $\qq$}\}$; then any extension of $w_p$ to
$L$ is a semihenselian valuation.
Note that if $v$ on $K$ is discrete, i.e., $\Gamma_K\cong \zz$,
then $w$ on $L$ cannot be semihenselian, since $\Gamma_L$ has no
nontrivial divisible subgroup; so, $w$ on $L$ must be henselian.
This preservation of the henselian property for discrete valuations
was asserted in \cite[Lemma, p.~195]{y}, but the proof given there is invalid.
\end{remark}

One associates to a valued division algebra $D$ a graded division algebra
as follows:
For each $\gamma\in \Gamma_D$, let
\begin{align*}
 D^{\ge\ga}  \ &=  \
\{ d \in D^{\ast} \mid v(d) \geq \ga \}\cup\{0\}, \text{ an additive
subgroup of $D$}; \qquad \qquad\qquad\qquad\qquad \ \\
D^{>\ga}  \ &=  \ \{ d \in D^{\ast} \mid v(d) > \ga \}\cup\{0\},
\text{ a subgroup
of $D^{\ge\ga}$};   \text{ and}\\
 \gr(D)_\gamma \ &= \
D^{\ge\ga}\big/D^{>\ga}.
\end{align*}
Then define
$$
 \gr(D)  \ = \  \textstyle\bigoplus\limits_{\ga \in \Ga_D}
\gr(D)_\gamma. \ \
$$
Because $D^{>\ga}D^{\ge\de} \,+\, D^{\ge\ga}D^{>\de}
\subseteq D^{>(\ga +
\de)}$ for all $\ga , \de \in \Ga_D$, the  multiplication on
$\gr(D)$ induced by multiplication on $D$ is
well-defined, giving that $\gr(D)$ is a graded  ring, called the
{\it associated graded ring} of $D$. The
multiplicative property
(1) of  the valuation $v$ implies that $\gr(D)$ is a graded
division ring.
Clearly,
we have ${\gr(D)}_0 = \overline{D}$ and $\Ga_{\gr(D)} = \Ga_D$.
For $d\in D^*$, we write $\widetilde d$ for the image
$d + D^{>v(d)}$ of $d$ in $\gr(D)_{v(d)}$.  Thus,
the map given by $d\mapsto \widetilde d$ is
a group epimorphism $\rho: D^* \rightarrow {\gr(D)^*}$ with
kernel~$1+M_D$, giving us the short exact sequence
\begin{equation}\label{grmap}
1 \ \longrightarrow 1+M_D  \ \longrightarrow D^* \  \longrightarrow  \ \gr(D)^*
 \ \longrightarrow  \ 1,
\end{equation} which will be used throughout. For a detailed
study of the associated graded algebra of a valued
division algebra
refer to \cite[\S4]{hwcor}.  As shown in \cite[Cor.~4.4]{hazwadsworth},
the  reduced norm maps
for $D$ and $\gr(D)$ are related by
\begin{equation}\label{nrdrel}
\widetilde{\Nrd_D(a)}  \ = \ \Nrd_{\gr(D)}(\widetilde a) \ \ \ \text{for all }
a\in D^*.
\end{equation}

Now let $K$ be a field with a henselian valuation $v$ and, as before,
let $D$ be a division algebra with center~$K$. Then $v$ extends uniquely
to a valuation on $D$, also denoted $v$,
 and one obtains associated to~$D$ the graded division algebra $\,\gr(D)=
\bigoplus_{\ga \in \Ga_D} D_{\ga}$. Further, suppose $D$ is
tame  with respect
to $v$. This implies that $[\gr(D):\gr(K)] = [D:K]$, $\gr(K) = Z(\gr(D))$ and
$D$ has a maximal subfield $L$ with $L$ tamely ramified over
$K$~(\cite[Prop.~4.3]{hwcor}).
We can then  associate to an
involution $\tau$ on $D$, a graded involution $\widetilde \tau$ on $\gr(D)$.
First, suppose
$\tau$ is of the first kind on $D$. Then $v \circ \tau$ is
also a valuation on $D$ which restricts to $v$ on $K$;
then, $v \circ \tau = v$
 since $v$ has a
unique extension to $D$. So,
$\tau$ induces a well-defined map $\widetilde{\tau}\colon \gr(D) \ra \gr(D)$,
defined on homogeneous elements by $\widetilde{\tau}(\widetilde{a}) =
\widetilde{\tau(a)}$ for all $a \in D^{\ast}$. Clearly,
$\widetilde{\tau}$ is a well-defined
graded involution on $\gr(D)$; it is of the first kind, as
it leaves $Z(\gr(D))=\gr(K)$ invariant.

If $\tau$ is a unitary involution on $D$, let $F=K^\tau$. In
this case, we need to assume  that the restriction of the valuation $v$
from $K$ to $F$
induces a henselian valuation on $F$, and that $K$ is tamely ramified
over $F$. Since $(v\circ\tau)|_F = v|_F$, an argument similar to the
one above
 shows that $v \circ \tau$ coincides with $v$ on $K$ and
thus on~$D$,  and the induced  map $\widetilde \tau$  on $\gr(D)$
as above is a graded
involution. That $K$ is tamely ramified over $F$
means that $[K:F] = [\gr(K):\gr(F)]$, $\overline{K}$ is separable
over $\overline{F}$, and $\chr(\overline F)
\nmid  |\Gamma_K : \Gamma_F|$.  Since $[K:F]=2$, $K$ is always tamely
ramified over $F$ if $\chr(\overline F) \not = 2$. But if
$\chr(\overline{F})=2$, $K$ is tamely ramified over $F$ if and
only if $[\overline{K}: \overline{F}]= 2$, $\Ga_K = \Ga_F$, and
$\overline{K}$ is separable (so Galois) over $\overline{F}$.
Since $K$ is Galois over $F$, the canonical map $\Gal(K/F)\to
\Gal(\ov K/\ov F)$ is surjective, by \cite[pp.~123--124, proof of
Lemma~5.2.6(1)]{EP}.
Hence,  $\tau$~induces the
nonidentity $\overline{F}$-automorphism $\ov \tau$ of $\overline{K}$.
Also $\widetilde \tau$ is unitary, i.e.,
$\widetilde{\tau}|_{\gr(K)}
\neq \id$. This is obvious if $\chr (\overline{F}) \neq 2$,
since then $K = F(\sqrt c)$ for some $c \in F^{\ast}$, and
$\widetilde{\tau}(\widetilde{\sqrt c}) =
\widetilde{\tau (\sqrt c)}=
-\widetilde{ \sqrt c} \neq \widetilde{\sqrt c}$.
If $\chr(\overline{F})= 2$, then $K$ is unramified over $F$ and
$\widetilde{\tau}|_{\gr(K)_0}= \ov{\tau}$ (the automorphism
of $\overline{K}$ induced by $\tau|_K$) which is nontrivial as
$\Gal(K/F)$ maps onto $\Gal(\overline{K}/ \overline{F})$; so
again $\widetilde{\tau}|_{\gr(K)} \neq \id$.
Thus, $\widetilde{\tau}$ is a  unitary graded involution
in any characteristic.
Moreover, for the graded fixed field $\gr(K)^{\widetilde{\tau}}$ we have
$\gr(F) \subseteq \gr(K)^{\widetilde{\tau}} \subsetneqq \gr(K)$ and
$[\gr(K):\gr(F)] =2$, so $\gr(K)^{\widetilde{\tau}} = \gr(F)$.

\begin{thm}\label{involthm1}
Let $(D,v)$ be a tame valued division algebra over a
henselian field $K$,
with $\chr(\overline K)\not = 2$.
 If  $\tau$ is an involution of the first kind on $D$,  then
$$
K_1(D,\tau)  \ \cong \  K_1(\gr(D),\widetilde \tau),
$$
and if  $\tau$ is symplectic, then
$$
K_1\Spin(D,\tau)  \ \cong  \ K_1\Spin(\gr(D),\widetilde \tau).
$$
\end{thm}

\begin{proof}

 Let $\rho\colon D^* \rightarrow \gr(D)^*$ be the
group epimorphism given in (\ref{grmap}). Clearly
$\rho(S_{\tau}(D))
\subseteq S_{\widetilde{\tau}}(\gr(D))$, so
$\rho (\Sigma_{\tau}(D))
\subseteq \Sigma_{\widetilde{\tau}}(\gr(D))$.
Consider the following diagram:

\begin{equation} \label{exactdiagramfirst}
\begin{split}
\xymatrix{
1 \ar[r] &  (1+M_D) \cap \Sigma_{\tau}(D)D' \ar[r] \ar[d]&
\Sigma_{\tau}(D)D' \ar[r]^-{\rho} \ar[d] &
\Sigma_{\widetilde{\tau}}(\gr(D))\gr(D)' \ar@{.>}[r] \ar[d]& 1  \\
1 \ar[r] & (1+M_D)  \ar[r] &
D^*  \ar[r]^-{\rho} & \gr(D)^*
 \ar[r] & 1. }
\end{split}
\end{equation}

The top row of the diagram is exact. To see this, note that
$\rho(D')=\gr(D)'$. Thus, it suffices to show that $\rho$ maps
$S_{\tau}(D) \cap D^{\ast}$ onto $S_{\widetilde{\tau}}(\gr(D)) \cap
\gr(D)^{\ast}$. For this, take any $d \in D^{\ast}$ with $\widetilde{d} =
\widetilde{\tau}(\widetilde{d})$.  Let ${b = \frac{1}{2} (d+ \tau(d)) \in
S_{\tau}(D)}$. Since $v(b) = v(\tau(b))$ and $\widetilde{d} +
\widetilde{\tau(d)} = 2 \widetilde{d} \neq 0$, ${\widetilde{b} = \frac{1}{2}
(\widetilde{d+ \tau(d)}) = \frac{1}{2} (\widetilde{d}+
\widetilde{\tau(d)}) = \widetilde{d}}$. Since $\tau$ on $D$ is an
involution of the first kind, the index of $D$ is a power of $2$
(\cite[Th.~1, \S 16]{draxl}). As ${\chr(\overline K)\not = 2}$, it follows
that the valuation is strongly tame, and by \cite[Lemma~2.1]{haz},
$$
{1+M_D \ = \ (1+M_K)[D^*,1+M_D]   \ \subseteq  \ \Sigma_{\tau}(D)D'}.
$$
Therefore,
the left vertical map is the identity map. It follows (for example
using the snake lemma) that $K_1(D,\tau) \cong K_1(\gr(D),\widetilde
\tau)$. The proof for  $K_1\Spin$  when $\tau$ is of symplectic type
is similar.
\end{proof}

The key to proving the corresponding result for unitary
involutions is the Congruence Theorem:

\begin{thm}[Congruence Theorem]\label{congthm}
Let $D$ be a tame  division algebra over a
field $K$ with henselian valuation $v$. Let
$D^{(1)} = \{ a\in  D^*\mid \Nrd_D(a) = 1\}$.
Then,
$$
D^{(1)} \, \cap  \, (1+M_D)  \ \subseteq  \ [D^*,D^*].
$$
\end{thm}

This theorem was proved by Platonov in \cite{platonov} for
$v$ a complete discrete valuation, and it was an essential tool
in all his calculations of $\SK$ for division rings.
The Congruence Theorem  was asserted by
Ershov in \cite{ershov}  in the generality given here.
A full proof is given in \cite[Th.~ B.1]{hazwadsworth}.

\begin{prop}[Unitary Congruence Theorem]\label{unitarycongthm}
Let $D$ be a tame  division algebra over a
field~$K$ with henselian valuation $v$, and let
 $\tau$ be a unitary involution  on $D$.  Let $F=K^\tau$.
If $F$ is henselian with respect to  $v|_F$
and $K$ is tamely ramified over $F$, then
$$
(1+M_D) \, \cap \,
\Sigma'_{\tau}(D)  \ \subseteq \  \Sigma_{\tau}(D).
$$
\end{prop}
\begin{proof}  The only published proof of this we know  is
\cite[Th.~4.9]{y}, which is just for the case $v$ discrete rank~1;
that proof is rather hard to follow, and appears to apply for other
 valuations only if $D$ is inertially split. Here we provide another
proof, in full generality.

We use the well-known facts that
\begin{equation}\label{nrd1+m}
\Nrd_D(1+M_D) \ = \ 1+M_K \text{\ \ \ and  \ \ \ }
N_{K/F}(1+M_K)  \ = \  1+M_F.
\end{equation}
(The second equation holds as $K$ is tamely ramified over $F$.)
See \cite[Prop.~2]{ershov} or \cite[Prop.~4.6, Cor.~4.7]{hazwadsworth}
for a proof.

Now, take $m \in M_D$ with $\Nrd_D(1 +m) \in F$. Then $\Nrd_D(1+m) \in F
\cap(1+M_K) = 1 +M_F$. By \eqref{nrd1+m} there is $c \in 1 +M_K$ with
$\Nrd_D(1 +m) = N_{K/F}(c) = c \tau(c)$, and there is $b \in 1+M_D$
with $\Nrd_D(b) = c$. Then,
$$
\Nrd_D(b \tau(b))  \ = \  c \tau(c)  \ = \  N_{K/F}(c)  \ = \  \Nrd_D(1+m).
$$
Let $s = (1+m)(b \tau(b))^{-1} \in 1+M_D$. Since $\Nrd_D(s) = 1$, by
the Congruence Theorem for $\SK$, Th.~\ref{congthm} above,
$s \in [D^{\ast}, D^{\ast}]
\subseteq \Sigma_{\tau}(D)$, (recall \eqref{primeinsigma})
. Since $b \tau(b) \in
S_{\tau}(D)$, we have $1+m = s(b \tau(b)) \in \Sigma_{\tau}(D)$.
\end{proof}

\begin{thm}\label{involthm2}
Let $D$ be a tame  division algebra over a
 field $K$ with henselian valuation $v$.
Let $\tau$ be a unitary involution  on $D$, and let $F=K^\tau$.
If $F$ is henselian with respect to $v|_F$
and $K$ is tamely ramified over $F$, then $\tau$ induces a
unitary graded involution $\widetilde{\tau}$ of $\gr(D)$  with
$\gr(F)=\gr(K)^{\widetilde{\tau}}$, and
$$
\SK(D, \tau)  \ \cong \  \SK(\gr(D), \widetilde{\tau}).
$$
\end{thm}

\begin{proof}
That
$\widetilde{\tau}$ is a  unitary graded involution
on $\gr(D)$ and $\gr(F)=\gr(K)^{\widetilde \tau}$  was
already observed (see the discussion before  Th.~\ref{involthm1}).
For the canonical
epimorphism $\rho\colon D^{\ast} \ra \gr(D)^{\ast}$,
$a \mapsto \widetilde{a}$, it follows from \eqref{nrdrel} that
$\rho( \Sigma'_{\tau}(D) \subseteq
\Sigma'_{\widetilde{\tau}}(\gr(D))$. Also, clearly $\rho(S_{\tau}(D))
\subseteq S_{\widetilde{\tau}}(\gr(D))$, so $\rho (\Sigma_{\tau}(D))
\subseteq \Sigma_{\widetilde{\tau}}(\gr(D))$. Thus, there is a
commutative diagram
\begin{equation} \label{exactdiagram}
\begin{split}
\xymatrix{1 \ar[r] & (1+M_D) \cap \Sigma_{\tau}(D) \ar[d] \ar[r] &
\Sigma_{\tau}(D) \ar[d] \ar[r]^-{\rho} &
\Sigma_{\widetilde{\tau}}(\gr(D))
\ar[d] \ar@{.>}[r] & 1 \\
1 \ar[r] &  (1+M_D) \cap \Sigma'_{\tau}(D) \ar[r] &
\Sigma'_{\tau}(D) \ar[r]^-{\rho} &
\Sigma'_{\widetilde{\tau}}(\gr(D)) \ar@{.>}[r] & 1,}
\end{split}
\end{equation}
where the vertical maps are inclusions, and the left vertical map is
bijective, by Prop.~\ref{unitarycongthm} above.

To see that the bottom row of  diagram~\eqref{exactdiagram} is exact at
$\Sigma'_{\widetilde{\tau}}(\gr(D))$, take $b \in D$ with
$\Nrd_{\gr(D)}(\widetilde{b}) \in \gr(F)$. Let $c = \Nrd_D(b) \in
K^{\ast}$. Then $\widetilde{c} = \Nrd_{\gr(D)} (\widetilde{b}) \in \gr(F)$, so
$\widetilde{c} = \widetilde{t}$ for some $t \in F^{\ast}$. Let $ u = c^{-1}
t \in 1 +M_K$. By \eqref{nrd1+m} above, there is $d \in 1+M_D$ with $\Nrd_D(d)
= u$. So, $\Nrd_D(bd) = cu = t \in F^{\ast}$. Thus, $bd \in
\Sigma'_{\tau}(D)$ and $\rho(bd) = \widetilde{\hsp bd\hsp } = \widetilde{b}$. This
gives the claimed exactness, and shows that the bottom row of
diagram~\eqref{exactdiagram} is exact.

To see that the top row of diagram~\eqref{exactdiagram} is exact at
$\Sigma_{\widetilde{\tau}}(\gr(D))$, it suffices to show that $\rho$ maps
${S_{\tau}(D) \cap D^{\ast}}$ onto $S_{\widetilde{\tau}}(\gr(D)) \cap
\gr(D)^{\ast}$. For this, take any $d \in D^{\ast}$ with $\widetilde{d} =
\widetilde{\tau}(\widetilde{d})$. If $\chr(\overline{F}) \neq 2$, as
in the proof of Th.~\ref{involthm1}, let $b =
\frac{1}{2} (d+ \tau(d)) \in S_{\tau}(D)$. Since $v(b) =
v(\tau(b))$ and $\widetilde{d} + \widetilde{\tau(d)} = 2 \widetilde{d} \neq
0$, we have ${\widetilde{b} = \frac{1}{2} (\widetilde{d+ \tau(d)}) = \frac{1}{2}
(\widetilde{d}+ \widetilde{\tau(d)}) = \widetilde{d}}$. If $\chr(\overline{F})
= 2$, then $K$ is unramified over $F$, so $\overline{K}$ is Galois over~
$\overline{F}$ with $[\overline{K}: \overline{F}]=2$, and the map
$\ov{\tau}\colon
\overline{K} \ra \overline{K}$ induced by $\tau$ is the nonidentity
$\overline{F}$-automorphism of $\overline{K}$. Of course,
$\overline{K} = \gr(K)_0$ and
$\ov{\tau} = \widetilde{\tau}|_{\gr(K)_0}$. Because $\overline{K}$
is separable
over $\overline{F}$, the trace $\textrm{tr}_{\overline{K}/\overline{F}}$ is
surjective, so
there is $r \in V_K$ with $\widetilde{r} + \widetilde{\tau}(\widetilde{r}) = 1
\in \gr(F)_0$. Let $c = rd + \tau(rd) \in S_{\tau}(D)$. We have
$\widetilde{\hsp rd\hsp } = \widetilde{r}\widetilde{d}$ and
$$
\widetilde{\tau(rd)} \, = \,
\widetilde{\tau}(\widetilde{\hsp rd\hsp})
 \, = \,  \widetilde{\tau}(\widetilde{r}\widetilde{d})  \, = \,
\widetilde{\tau}(\widetilde{d})\widetilde{\tau}(\widetilde{r})  \, = \,
\widetilde{\tau}(\widetilde{r})\widetilde{d}.
$$
 Since $v(rd) = v(\tau(rd))$ and
$\widetilde{\hsp rd\hsp } + \widetilde{\tau(rd)} = \widetilde{r}\widetilde{d}
+\widetilde{\tau}(\widetilde{r})\widetilde{d} = \widetilde{d} \neq 0$, we
have
$\widetilde{c} = \widetilde{\hsp rd\hsp } + \widetilde{\tau(rd)} = \widetilde{d}$. So,
in all cases $\rho(S_{\tau}(D) \cap D^{\ast}) = S_{\widetilde \tau}(\gr(D)) \cap
\gr(D)^{\ast}$, from which it follows that the bottom row of
diagram~\eqref{exactdiagram}
is exact. Since each row of~(\ref{exactdiagram}) is  exact,
we have a right exact sequence of
cokernels of the vertical maps, which yields the isomorphism of the
theorem.
\end{proof}

Having established the bridge between the unitary $K$-groups in
the graded setting and the non-graded henselian  case (Th.~\ref{involthm1},
Th.~\ref{involthm2}), we
can deduce
known formulas in the literature for the  unitary Whitehead
group of  certain valued division algebras,  by passing to
the graded setting.  The proofs are much easier than those
 previously available.
We will do this systematically for unitary  involutions
in Section~\ref{grsec}. Before we turn to that,  here is an example
with an involution of the first kind:

\begin{example}  Let $\sE$ be a graded division algebra over
its center $\sT$ with an involution $\tau$ of the first kind.
If $\sE$~is unramified over $\sT$, then, by using $\sE^*=\sE_0^*\sT^*$, it
follows easily  that
\begin{equation}\label{kuni}
K_1(\sE,\tau) \ \cong \  K_1(\sE_0,\tau|_{\sE_0}),
\end{equation} and, if $\chr(\sE) \ne 2$ and $\tau$ is symplectic,
\begin{equation}
K_1\Spin(\sE, \tau) \ \cong \  K_1\Spin(\sE_0,\tau|_{\sE_0}).
\end{equation}
Now if $D$ is a tame and unramified division algebra
over a henselian valued field and $D$ has an
 involution $\tau$ of the first kind, then the associated
graded division ring $\gr(D)$ is also unramified with the
corresponding graded involution $\widetilde \tau$ of the first kind;
then Th.~\ref{involthm1} and (\ref{kuni}) above show that
$$
K_1(D,\tau)  \ \cong  \ K_1(\gr(D),\widetilde{\tau})
 \ \cong \  K_1(\gr(D)_0,\tau|_{\gr(D)_0}) \ = \ K_1(\overline{D},
\overline{\tau}),
$$
yielding a theorem of Platonov-Yanchevski\u\i~\cite[Th. 5.11]{py85}
(that $K_1(D,\tau)\cong K_1(\overline{D},
\overline{\tau})$ when $D$ is unramified over $K$ and
the valuation is henselian and discrete rank $1$.)
Similarly, when $\chr(\ov D) \ne 2$ and $\tau$ is  symplectic,
$$
K_1\Spin(D,\tau)  \ \cong \  K_1\Spin(\gr(D),\widetilde{\tau})
 \ \cong \  K_1\Spin(\gr(D)_0,\tau|_{\gr(D)_0}) \ = \ K_1\Spin(\overline{D},
\overline{\tau}).
$$
\end{example}

\begin{remark}\label{goodrem}
We have the following commutative diagram connecting  unitary $\SK$ to
non-unitary $\SK$,
where $\SH(D,\tau)$ and $\SH(D)$ are the cokernels of $\Nrd\circ(1-\tau)$
and $\Nrd$ respectively (see diagram~\eqref{unnrd}).

\begin{equation}\label{goodd}
\begin{split}
\xymatrix{
1 \ar[r] & \SK(D,\tau) \ar[r] \ar[d] & D^*/\Sigma(D)
\ar[rr]^-{\Nrd\circ(1-\tau)} \ar[d]^{1-\tau} &
& K^* \ar[r] \ar[d]^{\id} & \SH(D,\tau) \ar[r]  \ar[d]& 1\\
1 \ar[r] & \SK(D) \ar[r]  &    D^*/D' \ar[rr]^-{\Nrd} &
& K^* \ar[r] & \SH(D) \ar[r] & 1.
}
\end{split}
\end{equation}

Now, let $D$ be a tame valued division algebra with  center $K$ and with
a unitary involution $\tau$, such that the valuation
restricts to a henselian valuation on $F=K^\tau$. By
Th.~\ref{involthm2}, $\SK(D,\tau)\cong \SK(\gr(D),\widetilde \tau)$
and by \cite[Th.~4.8, Th.~4.12]{hazwadsworth}, $\SK(D)\cong \SK(\gr(D))$
and $\SH(D)\cong\SH(\gr(D))$. However, $\SH(D,\tau)$ is not stable
under ``valued filtration'', i.e.,
$\SH(D,\tau)\not \cong \SH(\gr(D),\widetilde \tau)$. In fact
using~(\ref{grmap}), we can build a commutative diagram with exact rows,
\begin{equation*}
\begin{split}
\xymatrix{
1 \ar[r] & (1+M_K)\cap \Nrd(D^*)^{1-\tau} \ar[r] \ar@{^{(}->}[d] & \Nrd(D^*)^{1-\tau} \ar[r] \ar@{^{(}->}[d]& \Nrd\big(\gr(D)^*\big)^{1-\widetilde \tau} \ar[r]\ar@{^{(}->}[d] & 1\\
1 \ar[r] & 1+M_K \ar[r] & K^* \ar[r]  & \gr(K)^* \ar[r]  & 1,
}
\end{split}
\end{equation*}
which induces the exact sequence
$$
1\longrightarrow (1+M_K)\big/\big( (1+M_K)\cap\Nrd(D^*)^{1-\tau}\big)
\longrightarrow
\SH(D,\tau) \longrightarrow \SH(\gr(D),\widetilde \tau) \longrightarrow 1.
$$
By considering the norm $N_{K/F}\colon K^* \rightarrow F^*$, we clearly have
$\Nrd(D^*)^{1-\tau}
\subseteq \ker N_{K/F}$. However, by \eqref{nrd1+m},
${N_{K/F}\colon 1+M_K\rightarrow 1+M_F}$ is surjective, which shows that
$1+M_K\big/\big( (1+M_K)\cap\Nrd(D^*)^{1-\tau}\big)$ is not trivial and thus
${\SH(D,\tau)\not \cong \SH(\gr(D),\widetilde \tau)}$.
\end{remark}

\section{Graded Unitary $\SK$ Calculus}\label{grsec}

Let $\sE$ be a graded division algebra over its center $\sT$
with a unitary graded involution $\tau$,  and
let $\sR=\sT^\tau$.
Since $[\sT:\sR]=2=[\sT_0:\sR_0]\,|\Gamma_\sT:\Gamma_\sR|$,
there are  just two possible cases:

\begin{itemize}
\item[$\bullet$] $\sT$ is totally ramified over
$\sR$, i.e., $|\Gamma_\sT:\Gamma_\sR|=2$

\item [$\bullet$] $\sT$ is unramfied over
$\sR$, i.e., $|\Gamma_\sT:\Gamma_\sR|=1$.
\end{itemize}

We will consider $\SK(\sE,\tau)$ in these two cases separately
in \S\ref{totram} and \S\ref{unram}.

The following notation will be used throughout this section and the next:
Let $\tau'$ be another involution on $\sE$.  We write
$\tau' \sim \tau$ if
$\tau'|_{Z(\sE)}=\tau|_{Z(\sE)}$.  For $t \in \sE^*$,
let $\va_t$ denote  the map from $E$ to $E$ given by conjugation by $t$,
i.e., $\va_t(x)=txt^{-1}$.
Let $ \Sigma_0=\Sigma_\tau \cap \sE^*_0$ and
$\Sigma'_0=\Sigma_\tau'\cap \sE_0^*$.

We first  collect some facts which will be used below.
They all follow by easy calculations.

\begin{remarks}\hfill
\begin{enumerate}[\upshape(i)]\label{easyob}
\item \label{one} We have $\tau' \sim \tau$ if and only if there is a
$t\in \sE^*$ with $\tau(t)=t$ and $\tau'=\tau \va_t$.
(The proof is analogous to the ungraded version given, e.g. in
\cite[ Prop.~2.18]{kmrt}.)

\smallskip

\item \label{two} If $\tau'\sim \tau$, then $\Sigma_{\tau'}=\Sigma_{\tau}$
and $\Sigma'_{\tau'}=\Sigma'_{\tau}$; thus $\SK(\sE,\tau')=\SK(\sE,\tau)$.
(See~\cite[Lemma~1]{yin} for the analogous ungraded result.)
\smallskip
\item \label{three} For any $s \in \sE^*$, we have
$\tau \va_s=\va_{\tau(s)^{-1}}\tau$.  Hence, $\tau \va_s$ is an involution
(necessarily $\sim \tau$) if and only if $\tau \va_s=\va_{s^{-1}}\tau$
if and only if $\tau(s)/s \in \sT$.
\smallskip
\item \label{four} If $s \in \sE^*_\gamma$ and $\tau(s)=s$, then
$\Sigma_\tau' \cap \sE_\gamma= s \Sigma'_0$ and
$S_\tau\cap \sE_\gamma=s(S_{\tau_s}\cap \sE_0)$ where $\tau_s=\tau \va_s$.
\smallskip
\end{enumerate}
\end{remarks}

\subsection{$\sT/\sR$ totally ramified} \label{totram}
Let $\sE$ be a graded division algebra with a unitary graded involution $\tau$
 such that $\sT=Z(\sE)$ is totally
ramified over $\sR=\sT^\tau$. In this section we will show
that $\SK(\sE,\tau)=1$. Note that  the assumption that $\sT/\sR$ is
totally ramified implies that $\chr(\sT) \ne 2$.  For, if
$\chr(\sT) = 2$ and $\sT$ is totally ramified over a graded
subfield $\sR$ with $[\sT:\sR] = 2$, then for any $x\in \sT^*\setminus
\sR^*$, we have $\deg(x^2) \in \Gamma_R$, so $x^2\in \sR$;  thus, $\sT$ is
purely inseparable over $\sR$.  That cannot happen here,  as
$\tau|_\sT$ is a nontrivial $\sR$-automorphism of~$\sT$.

\begin{lemma}\label{six}
If $\sT$ is totally ramified over $\sR$, then
$\tau\sim\tau'$ for some graded involution $\tau'$, where
$\tau'|_{\sE_0}$  is of the first kind.
\end{lemma}

\begin{proof}
Let $Z_0=Z(\sE_0)$.
Since $\sT$ is totally ramified over $\sR$, $\sT_0=\sR_0$, so $\tau|_{Z_0}
\in \Gal(Z_0/\sT_0)$. Since the map
$\Theta_\sE\colon\Gamma_\sE \rightarrow \Gal(Z_0/\sT_0)$ is surjective
(see  \eqref{surj}), there is $ \gamma\in \Gamma_\sE$
with $\Theta_\sE(\gamma)=\tau|_{Z_0} $. Choose $y \in \sE^*_\gamma$ with
$\tau(y)=\pm y$. Then set $\tau'=\tau\va_{y^{-1}}$.
\end{proof}

\begin{example}\label{trex} Here is a construction of  examples of
 graded division algebras
$\sE$ with unitary  graded involution~$\tau$ with $E$ totally ramified over
$Z(\sE)^\tau$. We will see below that these are all such examples.
Let $\sR$ be any  graded field
with $\chr(\sR)\ne 2$, and let $\sA$ be a graded division algebra
with center $\sR$, such that $\sA$~is totally ramified over $\sR$ with
$\exp(\Gamma_\sA/\Gamma_\sR) = 2$.  Let $\sT$ be a graded field
extension of $\sR$ with $[\sT:\sR] = 2$, $\sT$ totally ramified over $\sR$,
and $\Gamma_\sT\cap \Gamma_\sA = \Gamma_\sR$.  Let $\sE = \sA\otimes_\sR
\sT$,  which is a graded central simple algebra over $\sT$, as
$\sA$ is graded central simple over $\sR$, by \cite[Prop.~1.1]{hwcor}.
But because $\Gamma_\sT\cap \Gamma_\sA = \Gamma_\sR$, we have
$\sE_0 = \sA_0 \otimes_{\sR_0} \sT_0 = \sR_0 \otimes_{\sR_0} \sR_0 = \sR_0$.
Since $\sE_0$ is a division ring, $\sE$ must be a graded division ring, which is
totally ramified over $\sR$, as $\sE_0 = \sR_0$.  Now, because $\sA$ is
totally ramified over $\sR$, we  have $\exp(\sA) = \exp(\Gamma_\sA/
\Gamma_\sR) = 2$, and $\sA= \sQ_1\otimes_\sR \ldots \otimes_\sR\sQ_m$, where
each $\sQ_i$ is a  graded symbol algebra of degree at most $2$, i.e., a
graded quaternion algebra.  Let $\sigma_i$ be a graded involution of the
first kind on $\sQ_i$ (e.g., the canonical symplectic graded involution),
 and let $\rho$ be the nonidentity $\sR$-automorphism of $\sT$.
Then, $\sigma = \sigma_1\otimes \ldots \otimes \sigma_m$ is a graded involution
of the first  kind on $\sA$, so $\sigma\otimes \rho$  is a unitary graded
involution on~$\sE$, with $\sT^\tau = \sR$.

\end{example}

\begin{prop}\label{total}
If $\sE$  is totally ramified over $\sR$, and $\sE \ne  \sT$,
then $\Sg=\sE^*$, so
$\SK(\sE,\tau)=1$.  Furthermore, $\sE$ and $\tau$ are as described
in Ex.~\ref{trex}.
\end{prop}

\begin{proof}
We have $\sE_0=\sT_0=\sR_0$.
For any  $\gamma \in \Gamma_\sE$, there is a nonzero $a\in \sE_\gamma$
with $\tau(a) = \epsilon a$ where $\epsilon = \pm 1$.  Then,
for any $b\in \sE_\gamma$, $b = ra$ for some $r\in  \sE_0 = \sR_0$.
Since $r$ is  central and symmetric, $\tau(b) = \epsilon b$.
Thus,  every element of $\sE^*$ is symmetric or
skew-symmetric. Indeed, fix any $t\in \sT^*\setminus \sR^*$.
Then $\tau(t) \ne t$, as $t\notin \sR^*$.  Hence, $\tau(t) =
-t$.  Since $t$ is central and skew-symmetric, every
$a\in \sE^* $ is symmetric  iff $ta$ is  skew-symmetric.
Thus,  $\sE^* = S_\tau^* \cup tS_\tau^*$.  To see that
$\Sigma_\tau = \sE^*$, it suffices to show that $t\in  \Sigma_\tau$.
To see this, take any $c,d\in \sE^*$ with $dc\ne cd$.  (They exist, as
$\sE \ne \sT$.)  By replacing $c$ (resp.~$d$) if necessary by
$tc$ (resp.~$td$), we may assume that $\tau(c) = c$ and $\tau(d) = d$.
Then, $dc = \tau(cd) = \epsilon cd$, where $\epsilon = \pm 1$;
since $dc\ne cd$, $\epsilon = -1$;  hence $\tau(tcd) = tcd$.
Thus, $t = (tcd) c^{-1}  d^{-1} \in \Sigma_\tau(\sE)$, completing the proof
that $\Sigma_\tau(\sE)  = \sE^*$.

For $\gamma \in \Gamma_\sE$, let $\ov \gamma = \gamma+\Gamma_\sT \in
\Gamma_\sE/\Gamma_\sT$.
To see the structure of $\sE$, recall that as $\sE$ is totally ramified
over~$\sT$ there is a well-defined nondegenerate $\zz$-bilinear symplectic
pairing $\beta\colon (\Gamma_\sE/\Gamma_\sT) \times\Gamma_\sE/\Gamma_\sT)
\to \sE_0^*$ given by $\beta(\ov \gamma, \ov \delta) = y_\gamma y_\delta
y_\gamma^{-1}y_\delta^{-1}$ for any nonzero $y_\gamma\in  \sE_\gamma$,
$y_\delta\in \sE_\delta$.  The computation above for $c$ and $d$ shows
that $\im (\beta) = \{\pm 1\}$.  Since the pairing $\beta$ is nondegenerate
by \cite[Prop.~2.1]{hwcor}  there is a symplectic base  of $\Gamma_\sE/\Gamma_\sT$, i.e., a subset
$\{\ov \gamma_1, \ov \delta_1, \ldots , \ov \gamma_m, \ov \delta_m\}$   of
$\Gamma_\sE/\Gamma_\sT$ such that $\beta(\ov \gamma_i, \ov \delta_i) = -1$
while $\beta(\ov \gamma_i, \ov \gamma_j)  = \beta(\ov \delta_i, \ov \delta_j)
=1$ for all $i,j$, and $\beta(\ov \gamma_i, \ov \delta_j) = 1$ whenever
$i\ne j$, and $\Gamma_\sE = \langle \gamma_1, \delta_1,\ldots,
\gamma_m, \delta_m\rangle + \Gamma_\sT$.  Choose any nonzero $\bi_i\in
\sE_{\gamma_i}$ and $\bj_i\in \sE_{\delta_i}$.  The properties of the
$\ov\gamma_i, \ov \delta_i$  under $\beta$ translate  to:
$\bi_i \bj_i = - \bj_i \bi_i$ while $\bi_i \bi_j =  \bi_j \bi_i$
and $\bj_i \bj_j =  \bj_j \bj_i$ for all $i,j$, and
$\bi_i \bj_j =  \bj_j \bi_i$ whenever $i \ne j$. Since $\beta(
\ov {2\gamma_i}, \ov \eta) = 1$ for all $i$ and all $\eta\in \Gamma_\sE$,
each~$\bi_i^2$ is central in $\sE$.  But also $\tau(\bi_i^2) = \bi_i^2$,
as $\tau(\bi_i) = \pm \bi_i$. So, each $\bi_i^2 \in \sR^*$, and likewise each
$\bj_i^2\in \sR^*$.  Let
$\sQ_i = \sR\text{-span}\{1, \bi_i,  \bj_i,  \bi_i\bj_i\}$ in~$\sE$.
The relations on the $\bi_i, \bj_i$ show that each $\sQ_i$ is a
graded quaternion algebra over $\sR$, and the distinct~$\sQ_i$ centralize
each other in $\sE$.    Since each $\sQ_i$ is graded central simple over $\sR$,
$\sQ_1\otimes_\sR \ldots \otimes_\sR\sQ_m$ is graded  central  simple over
$\sR$
by \cite[Prop.~1.1]{hwcor}.  Let $\sA = \sQ_1\ldots \sQ_m \subseteq \sE$.  The graded
$\sR$-algebra epimorphism $\sQ_1\otimes_\sR \ldots \otimes_\sR\sQ_m \to
\sA$ must be an isomorphism, as the domain is graded simple.
If $\Gamma_\sT \subseteq \Gamma_\sA$, then $\sT \subseteq \sA$,
since $\sE$ is totally ramified over $\sR$.  But this cannot occur, as
$\sT$ centralizes $\sA$ but $\sT\supsetneqq \sR = Z(\sA)$.  Hence,
as $|\Gamma_\sT:\Gamma_\sR| = 2$, we must have $\Gamma_\sT \cap \Gamma_\sA
 = \Gamma_\sR$.
The graded $\sR$-algebra homomorphism $\sA \otimes_\sR \sT \to
\sE$ is injective since its domain is graded simple, by \cite[Prop.~1.1]{hwcor};
it is also surjective, since $\sE_0 = \sR_0 \subseteq \sA\otimes_\sR \sT$
and $\Gamma_{\sA \otimes_\sR \sT} \supseteq \langle \gamma_1, \delta_1,\ldots,
\gamma_m, \delta_m\rangle +\Gamma_\sT = \Gamma_\sE$.  Clearly,
$\tau = \tau|_\sA \otimes \tau|_\sT$.
\end{proof}

\begin{prop}\label{completely}
If $\sE \ne \sT$ and $\sT$  is totally ramified over $\sR$, then
$\Sg=\sE^*$, so $\SK(\sE,\tau)=1$.
\end{prop}

\begin{proof}
The case where $\sE_0 = \sT_0$ was covered by Prop.~\ref{total}.
Thus, we may assume that $\sE_0 \supsetneqq\sT_0$.
By Lemma~\ref{six} and Remark~\ref{easyob}(\ref{two}), we can assume
that $\tau|_{\sE_0}$ is of the first kind. Further, we can assume that
$\sE_0^*=\Sigma_{\tau|_{\sE_0}}(\sE_0).$
For, if $\tau|_{\sE_0}$
is symplectic, take any $a\in \sE_0^*$ with $\tau(a)=-a$, and let
$\tau'=\tau\va_a$.
Then, $\tau' \sim \tau$ (see Remark~\ref{easyob}(\ref{three})).
 Also, $\tau'|_{Z(\sE_0)}=\tau|_{Z(\sE_0)}$,  as $a \in \sE_0$ and so
$\va_a |_{Z(\sE_0)}=\id$. Therefore, $\tau'|_{\sE_0}$ is of the first kind.
But as $\tau(a)=-a$, $\tau'|_{\sE_0}$ is orthogonal.
Thus $\sE_0^*=\Sigma_{\tau'|_{\sE_0}}(\sE_0)$, as noted at the beginning of
\S\ref{firstk}. Now replace $\tau$ by $\tau'$.

We consider two cases.

{Case I.} Suppose for each $\gamma\in \Gamma_\sE$ there is
$x_\gamma \in \sE^*_\gamma$ such that $\tau(x_\gamma)=x_\gamma$.
 Then, $\sE^* = \bigcup_{\gamma \in \Gamma_\sE}\sE_0^*x_\gamma
\subseteq \Sigma_\tau(\sE)$, as desired.

{Case II.} Suppose there is $\gamma\in \Gamma_\sE$ with
$\sE_\gamma \cap S_\tau=0$. Then $\tau(d)=-d$ for each $d\in \sE_\gamma$.
Fix $t \in \sE^*_\gamma$. For any $a \in \sE_0$, we have $ta \in \sE_\gamma$;
 so,
$-ta=\tau(ta)=\tau(a)\tau(t)=-\tau(a)t$.  That is,
\begin{equation}\label{localeq}
\tau(a) \, = \, \va_t(a)\quad\text{for all }a \in \sE_0.
\end{equation}
Let $\tau''=\tau\va_t$, which is a unitary involution on $\sE$ with
$\tau''\sim\tau$ (see Remark~\ref{easyob}(\ref{three})). But,
$\tau''(a)=a$ for all
$a \in \sE_0$, i.e., $\tau''|_{\sE_0}=\id$. This implies that $\sE_0$ is
a field. Replace $\tau$ by $\tau''$. The rest of the argument uses this
new $\tau$.
So $\tau|_{\sE_0}=\id$. If we are now in
Case I for this $\tau$, then we are done by Case I. So, assume we are
in Case II.
Take any $\gamma\in \Gamma_\sE$ with
$\sE_\gamma \cap S_\tau=0$.  For any nonzero $t\in \sE_\gamma$,
equation~(\ref{localeq}) applies to $t$, showing $\va_t(a)=\tau(a)=a$
for all $a \in \sE_0$; hence for the map $\Theta_\sE$ of
\eqref{surj}, $\Theta_\sE(\gamma) = \id_{\sE_0}$. But recall that
$\sE_0$~is Galois over $\sT_0$ and $\Theta_\sE\colon \Gamma_\sE\to
\Gal(\sE_0/\sT_0)$ is surjective.  Since $\sE_0 \ne \sT_0$,
  there is $\delta \in \Gamma_\sE$ with
$\Theta_\sE(\delta) \ne \id$.  Hence, there must be some
$s\in \sE_\delta^*\cap S_\tau$.  Likewise, since $\Theta_\sE(\gamma-\delta)
= \Theta_\sE(\gamma) \Theta_\sE(\delta)^{-1} \ne \id$, there
is some $r\in \sE_{\gamma-\delta}^*\cap S_\tau$.  Then,
as $rs\in \sE_\gamma^*$, we have $\sE_\gamma^*  = \sE_0^*rs \subseteq
\Sigma_\tau$.  This is true for every $\gamma$  with $\sE_\gamma\cap S_\tau
= 0$.    But for any other $\gamma\in \Gamma_\sE$, there is an $x_\gamma$
in  $\sE_\gamma^* \cap S_\tau$; then $\sE_\gamma^* = \sE_0^*x_\gamma \subseteq
\Sigma_\tau$. Thus, $\sE^* = \bigcup_{\gamma\in \Gamma_\sE}\sE_\gamma^*
\subseteq \Sigma_\tau$.
 \end{proof}

\subsection{$\sT/\sR$    unramified} \label{unram}

Let $\sE$ be a graded division algebra with a unitary involution
 $\tau$ such that $\sT=Z(\sE)$ is unramified over $\sR=\sT^\tau$.
In this subsection, we will give a general formula for $\SK(\sE,\tau)$ in
terms of data in $\sE_0$.

\begin{lemma}\label{unramfacts}
Suppose $\sT$ is unramified over $\sR$.  Then,
\begin{enumerate}[\upshape (i)]
 \item \label{seven}
Every $\sE_\gamma$ contains both nonzero symmetric
and skew symmetric elements.

\smallskip

  \item  \label{eight}
$Z(\sE_0)$ is a generalized dihedral extension for
$\sT_0$ over
$\sR_0$ $($see Def.~\ref{gendi}$)$.

\smallskip

  \item  \label{nine}
If $\sT$ is unramified over $\sR$, then
$\SK(\sE,\tau)=\Sigma'_0/\Sigma_0$.
\end{enumerate}
\end{lemma}

\begin{proof}\hfill

(i) If $\chr(\sE)=2$, it is easy to see that every
$\sE_\gamma$ contains a symmetric element (which is also  skew
symmetric)  regardless of any assumption on $\sT/\sR$.
Let $\chr(\sE)\not=2$.
Since $[\sT_0:\sR_0]=2$ and $\sR_0=\sT_0^\tau$,  there is
$ c\in \sT_0$ with $\tau(c)=-c$. Now there is
$t\in \sE_\gamma$, $t\not =0$, with $\tau(t)=\epsilon t$ where
$\epsilon=\pm 1$. Then $\tau(c t)=-\epsilon c t $.

(ii) Let $G=\Gal(Z(\sE_0)/\sR_0)$ and $H=\Gal(Z(\sE_0)/\sT_0)$.
Note that $[G:H]=2$. Since $\tau$ is unitary,
$\tau |_{Z(\sE_0)} \in G\setminus H$. We will denote
$\tau |_{Z(\sE_0)}$ by $\overline \tau$ and will show that for
any $h\in H$, $(\overline \tau h)^2=1$. By (\ref{surj}),
$\Theta_\sE\colon \Gamma_\sE \rightarrow \Gal(Z(\sE_0)/\sT_0)$ is onto,
so there is $\gamma \in \Gamma_\sE$,  such that $\Theta_\sE(\gamma)=h$.
Also by (\ref{seven}),
there is an $x\in \sE_\gamma^*$ with $\tau(x)=x$.
Then $\tau\va_x$ is an involution,  where $\va_x$ is  conjugation by $x$;
therefore, $\tau \va_x|_{Z(\sE_0)} \in G$ has order $2$.
But $\va_x|_{Z(\sE_0)}=\Theta_\sE(\gamma)=h$. Thus
$(\overline \tau h)^2=1$.

(iii) By~(\ref{seven}), for each $\gamma \in \Gamma_\sE$, there is
$s_\gamma \in \sE_\gamma$,  $s_\gamma\not =0$, with
$\tau(s_\gamma)=s_\gamma$. By Remark~\ref{easyob}\eqref{four},
${\Sigma_\tau'=\bigcup_{\gamma \in \Gamma_\sE} s_\gamma \Sigma'_0}$.
Since each $s_\gamma \in S_\tau \subseteq \Sigma_\tau$, the injective map
$\Sigma'_0/\Sigma_0 \rightarrow \Sigma_\tau'/\Sigma_\tau$ is an isomorphism.
\end{proof}

To simplify notation in the next theorem, let $\ov \tau =
\tau|_{Z(\sE_0)} \in \Gal(Z(\sE_0)/\sR_0)$, and
for any $h \in \Gal(Z(\sE_0)/\sT_0)$,
 write $\Sigma_{h \overline \tau}(\sE_0)$ for $\Sigma_\rho(\sE_0)$ for any unitary
involution $\rho$ on $\sE_0$ such that $\rho|_{Z(\sE_0)}=h\overline \tau$.
This is well-defined, independent of the choice of $\rho$, by the ungraded
analogue of Remark\eqref{easyob}\eqref{two}.

\begin{theorem}\label{msem}
Let $\sE$ be a graded division algebra with center $\sT$, with a unitary
graded involution $\tau$, such that $\sT$~is unramified over $\sR=\sT^\tau$.
For each $\gamma \in \Gamma_\sE$ choose a nonzero
$x_\gamma \in S_\tau \cap \sE_\gamma$. Let $H=\Gal(Z(\sE_0)/\sT_0)$.
Then,
$$
\SK(\sE,\tau) \ \cong  \ (\Sigma'_\tau \cap \sE_0) \big /
(\Sigma_\tau \cap \sE_0),
$$
with
\begin{equation}\label{ersh}
\Sigma'_\tau \cap \sE_0  \ = \
\big \{ a \in \sE_0^*\mid
N_{Z(\sE_0)/\sT_0}\Nrd_{\sE_0}(a)^\dlambda\in \sR_0 \big \},
\textrm{  \ \ \ where \ \ \ } \dlambda \, = \
\ind(\sE)/\big(\ind(\sE_0) \, [Z(\sE_0):\sT_0]\big)
\end{equation}
and
\begin{equation}\label{ersh1}
\Sigma_\tau\cap\sE_0 \ = \ P\cdot X, \textrm{\ \ \ where\ \ \ }
P\, =\ \textstyle{\prod}_{h \in H}\Sigma_{h\overline \tau}(\sE_0)
\textrm{ \ \ and \ \ }
X\, = \ \langle x_\gamma x_\delta x_{\gamma+\delta}^{-1} \mid
\gamma,\delta \in \Gamma_\sE \rangle \, \subseteq \, \sE_0^*.
\end{equation}
Furthermore, if $H=\langle h_1,\dots,h_m \rangle$, then
$P \, = \ \prod_{(\varepsilon_1,\dots,\varepsilon_m)\in \{0,1\}^m}
\Sigma_{h_1^{\varepsilon_1}\dots h_m^{\varepsilon_m} \overline \tau}(\sE_0)$.
\end{theorem}

Before proving the theorem, we record the following:

\begin{lemma}\label{lem5}
Let $A$ be a central simple algebra over a field $K$, with an
involution $\tau$ and an automorphism or anti-automorphism $\sigma$.
Then,

\begin{enumerate}[\upshape(i)]
\item  \label{hyh1} $\sigma \tau \sigma^{-1}$ is an involution of $A$ of the same
 kind as $\tau$, and
$$
S_{\sigma \tau \sigma^{-1}}  \, = \, \sigma(S_\tau),
\textrm{ \ \  so \  \  } \Sigma_{\sigma \tau \sigma^{-1}} \, = \,
\sigma(\Sigma_\tau).
$$

\item \label{hyh2} Suppose $A$ is a division ring.
If $\sigma$ and $\tau$ are each unitary involutions, then
$($writing {${S_\tau^*=S_\tau\cap A^*}$}$)$,
$$
S_\tau^* \ \subseteq  \ S_\sigma^* \cdot \sigma(S_\tau^*)
 \ = \ S_\sigma^* \cdot S_{\sigma \tau \sigma^{-1}}^*, \textrm{  \  \ so \  \ ~~}
\Sigma_\tau  \ \subseteq  \ \Sigma_\sigma\cdot \Sigma_{\sigma \tau \sigma^{-1}}.
$$
\end{enumerate}
\end{lemma}
\begin{proof}\hfill

(i) This follows by easy calculations.

(ii) Observe that if $a \in S_\tau^*$,
then $a=\big (a\sigma(a)\big)\sigma(a^{-1})$ with ${a\sigma(a) \in S_\sigma^*}$
and $\sigma(a^{-1}) \in \sigma(S_\tau^*)=S_{\sigma \tau \sigma^{-1}}^*$ by~(\ref{hyh1}).
Thus,~(\ref{hyh2}) follows from~(\ref{hyh1}) and the fact
that ${A' \subseteq \Sigma_\tau \cap \Sigma_\sigma}$
(see \eqref{primeinsigma}).
\end{proof}

\begin{proof}[Proof of Theorem~\ref{msem}]

First note that by Lemma~\ref{unramfacts}\eqref{nine} the canonical map
$$
(\Sigma'_\tau \cap \sE_0) \, \big/ \, (\Sigma_\tau \cap \sE_0)
 \ \longrightarrow  \ \Sigma'_\tau/\Sigma_\tau \ = \ \SK(\sE,\tau)
$$
is an isomorphism.
The description of $\Sigma'_\tau \cap \sE_0$
in~(\ref{ersh}) is immediate from the fact that
for $a \in \sE_0$,
$\Nrd_{\sE}(a)=N_{Z(\sE_0)/\sT_0}\Nrd_{\sE_0}(a)^\dlambda\in \sT_0$
(see Remark~\ref{grfacts}\eqref{rnrd}).

For $\Sigma_\tau \cap \sE_0$, note that for each
$\gamma \in \Gamma_\sE$, if $a\in \sE_0$, then
$a x_\gamma \in S_\tau$ if and only if
 $x_\gamma \tau(a) x_\gamma^{-1}=a$. That is,
$S_\tau \cap \sE_\gamma = S(\va_{x_\gamma} \tau;\sE_0)x_\gamma$,
where $S(\va_{x_\gamma} \tau;\sE_0)$ denotes the set of
symmetric elements in $\sE_0$ for the unitary involution
$\va_{x_\gamma} \tau|_{\sE_0}$.
Therefore,
$$
\Sigma_\tau \cap \sE_0  \ = \
\big \langle S(\va_{x_\gamma} \tau;\sE_0)^* x_\gamma \mid
\gamma \in \Gamma_\sE \big \rangle \, \cap \,  \sE_0.
$$
Take a product $a_1x_1\ldots a_kx_k$ in $\Sigma_\tau \cap \sE_0$
where each $x_i=x_{\gamma_i}$ for some $\gamma_i \in \Gamma_\sE$
and $a_i \in S(\va_{x_i} \tau;\sE_0)^*$. Then,
\begin{equation}\label{pet}
a_1x_1\dots a_k x_k \ = \
a_1\va_{x_1}(a_2)\ldots\va_{x_1\ldots x_{i-1}}(a_i)\ldots
\va_{x_1\ldots x_{k-1}}(a_k)x_1\ldots x_k \, \in \,
\sE_{\gamma_1+\ldots+\gamma_k}.
\end{equation}
So, $\gamma_1+\ldots+\gamma_k=0$. Now, as $a_i \in S(\va_{x_i} \tau;\sE_0)$
and $\tau \va_{x_j}^{-1}=\va_{x_j}\tau$ for all $j$,
by Lemma~\ref{lem5}(\ref{hyh1}) we obtain
\begin{equation}\label{longd}
\va_{x_1\ldots x_{i-1}}(a_i)\, \in \,
S(\va_{x_1}\ldots \va_{x_{i-1}}(\va_{x_i}\tau)
\va_{x_{i-1}}^{-1}\ldots\va_{x_1}^{-1};
\sE_0)^* \ = \ S(\va_{x_1\ldots x_{i-1}x_ix_{i-1}\ldots x_1} \tau;\sE_0)^*
 \ \subseteq  \ \Sigma_{h\overline \tau}(\sE_0)  \ \subseteq  \, P,
\end{equation}
where $h=\va_{x_1\dots x_{i-1}x_ix_{i-1}\dots x_1}|_{Z(\sE_0)} \in H$.
Note also that if $k=1$, then
$x_1 \in S_\tau \cap \sE_0^* \subseteq \Sigma_{\overline \tau}(\sE_0)
\subseteq P.$

If $k>1$, then
$$
x_1\ldots x_k \ = \ x_{\gamma_1}\ldots x_{\gamma_k} \
= \ (x_{\gamma_1} x_{\gamma_2}
x_{\gamma_1+\gamma_2}^{-1})(x_{\gamma_1+\gamma_2}x_{\gamma_3}
\ldots x_{\gamma_k}),
$$
with
$(\gamma_1+\gamma_2)+\gamma_3+\ldots+\gamma_k=0$. It follows by induction
on $k$ that $x_1 \dots x_k \in X$. With this and ~(\ref{pet})
and ~(\ref{longd}), we have $a_1x_1\dots a_kx_k \in P\cdot X$
(which is a group, as $\sE_0' \subseteq \Sigma_{\overline \tau}(\sE_0)
\subseteq  P$
by
\eqref{primeinsigma}), showing that $\Sigma_\tau \cap \sE_0
\subseteq P \cdot X$. For the reverse inclusion, take any
$h \in H$ and choose $\gamma \in \Gamma_\sE$ with
$\va_{x_\gamma}|_{Z(\sE_0)}=h$. Then,  $x_\gamma \in S_\tau^*
\subseteq \Sigma_\tau$ and $S(\va_{x_\gamma}\tau;\sE_0)^*x_\gamma=
S_\tau^*\cap \sE_\gamma \subseteq \Sigma_\tau$, so
$\Sigma_{h\overline \tau}(\sE_0)=\Sigma_{\va_{x_\gamma}\tau}(\sE_0)=
\langle S(\va_{x_\gamma}\tau;\sE_0)^* \rangle \subseteq
\Sigma_\tau \cap \sE_0$. Thus, $P \subseteq \Sigma_\tau \cap \sE_0$,
and clearly also $X \subseteq \Sigma_\tau \cap \sE_0$. Hence,
$\Sigma_\tau \cap \sE_0 =P\cdot X$.

The final equality for $P$ in the Theorem follows from
Lemma~\ref{lembe} below by taking $U=\sE_0^*$, $A=H$, and
$W_h=\Sigma_{h\overline \tau}(\sE_0)$ for $h \in H$. To see that
the lemma applies, note that each $\Sigma_{h\overline \tau}(\sE_0)$
contains $\sE_0'$ by \eqref{primeinsigma}.  Furthermore, take any $h,\ell \in H$, and choose
$x,y \in E^* \cap S_\tau$ with $\va_x|_{Z(\sE_0)}=h$ and
$\va_y|_{Z(\sE_0)}=\ell$. Then,
$$
(\va_y\tau)(\va_x\tau)(\va_y\tau)^{-1} \ = \
\va_y\tau\va_x\va_y^{-1} \ = \ \va_{yx^{-1}y}\tau,
$$
and $\va_{yx^{-1}y}|_{Z(\sE_0)}=\ell h^{-1}\ell=\ell^2h^{-1}$. Hence, by
Lemma~\ref{lem5}(\ref{hyh2}),
$\Sigma_{h\overline \tau}(\sE_0)\subseteq
\Sigma_{\ell\overline \tau}(\sE_0)\Sigma_{\ell^2h^{-1}\overline \tau}(\sE_0)$.
This shows that hypothesis~(\ref{hyp}) of Lemma~\ref{lembe} below is
satisfied here.
\end{proof}

\begin{lemma}\label{lembe}
Let $U$ be a group, $A$ an abelian group, and $\{W_a\mid a\in A\}$ a family
of subgroups of $U$ with each $W_a \supseteq [U,U]$. Suppose
\begin{equation}\label{hyp}
W_a \, \subseteq  \, W_bW_{2b-a} \textrm{ \ \ for~all  } a,b \in A.
\end{equation}
If $A=\langle a_1,\dots,a_m \rangle$, then
$$
\textstyle{\prod\limits_{a\in A}} W_a \ = \
\textstyle{\prod\limits_{(\varepsilon_1,\dots,\varepsilon_m)\in \{0,1\}^m}}
W_{\varepsilon_1a_1+\ldots+\varepsilon_ma_m}. $$
\end{lemma}
\begin{proof}
Since each $W_a \supseteq [U,U]$, we have $W_aW_b=W_bW_a$, and this is a
subgroup of $U$, for all $ a,b \in A$.
Let
$$
Q \, = \textstyle{\prod\limits_{(\varepsilon_1,\dots,\varepsilon_m)\in \{0,1\}^m}}
 W_{\varepsilon_1a_1+\ldots+\varepsilon_ma_m}.
$$
We prove by induction on $m$ that each $W_a\subseteq Q$. The lemma then
follows, as $Q$ is a subgroup of $U$.
Note that condition~(\ref{hyp}) can be conveniently restated,
\begin{equation}\label{hyp1}
\textrm{ if  \ } a+b \, = \, 2d \in A, \textrm{  \ then  \ }
W_a  \, \subseteq \,  W_dW_b.
\end{equation}
Take any $c\in A$. Then, (\ref{hyp1}) shows that $W_{-c}\subseteq W_0W_c$.
Take any $i \in \mathbb Z$, and suppose $W_{ic}\subseteq W_0W_c$. Then by
(\ref{hyp1}) $W_{-ic}\subseteq W_0W_{ic} \subseteq W_0W_c$. So, by
(\ref{hyp1}) again, $W_{(i+2)c}\subseteq W_cW_{-ic} \subseteq W_0W_c$ and
${W_{(i-2)c}\subseteq W_{-c}W_{-ic} \subseteq W_0W_c}$. Hence, by induction
(starting with $j=0$ and $j=1$), $W_{jc}\subseteq W_0W_c$ for every
$j \in \mathbb Z$. This proves the lemma when $m=1$.

Now assume $m>1$ and let $B=\langle a_1,\dots,a_{m-1} \rangle \subseteq A$.
By induction, for all
$b\in B$,
$$
W_b \ \subseteq
\textstyle{\prod\limits_{(\varepsilon_1,\dots,\varepsilon_{m-1})\in \{0,1\}^{m-1}}}
W_{\varepsilon_1a_1+\ldots+\varepsilon_{m-1}a_{m-1}} \ \subseteq  \, Q.
$$
Also, by the cyclic case done above,
$W_{j a_m} \subseteq W_0W_{a_m} \subseteq Q$
for all $j \in \mathbb Z$. So, for any $b \in B, j \in \mathbb Z$,
using~(\ref{hyp1}),
\begin{equation}\label{thte}
W_{2b+ja_m}  \, \subseteq  \, W_b W_{-ja_m}  \, \subseteq  \, Q,
\end{equation}
and
\begin{equation}\label{thte2}
W_{b+2ja_m}  \, \subseteq  \, W_{ja_m} W_{-b} \, \subseteq  \, Q.
\end{equation}
Let $d=a_{i_1}+\ldots+a_{i_\ell}$ for any indices
$1\leq i_1< i_2 <\ldots <i_\ell \leq m-1$. Since
$W_{d+a_m} \subseteq Q$ by hypothesis,
from~(\ref{hyp1}) and~(\ref{thte2}) it follows that
\begin{equation}\label{d+3}
W_{d+3a_m} \,  \subseteq  \ W_{d+2a_m} W_{d+a_m}  \, \subseteq  \, Q.
\end{equation}
Now, take any element of $A$; it has the form $b+ja_m$ for some
$b\in B$ and $j \in \mathbb Z$. If $b \in 2B$ or if $j$ is even,
then~(\ref{thte}) and~(\ref{thte2}) show that
$W_{b+ja_m} \subseteq Q$. The remaining case is that $j$ is
odd and $b \not \in 2B$, so $b=2c+d$, where $c\in B$ and
$d=a_{i_1}+\ldots+a_{i_\ell}$ for some indices with
$1\leq i_1<i_2<\ldots <i_\ell \leq m-1$. Set $q=1$ if
$j \equiv 3 \pmod{4}$ and $q=3$ if
$j \equiv 1  \pmod{4}$. Then, $W_{d+qa_m}\subseteq Q$ by definition
if $q=1$ or by \eqref{d+3} if $q=3$.
Hence, by~(\ref{hyp1}),
\begin{equation}\label{ww}
W_{b+ja_m} \, = \ W_{(2c+d)+ja_m}  \, \subseteq  \
W_{(c+d)+((j+q)/2)a_m}W_{d+qa_m} \,  \subseteq  \, Q,
\end{equation}
using~(\ref{thte2}) as $c+d \in B$ and $(j+q)/2$ is even.
Thus, $W_a \subseteq Q$, for all $a \in A$.
\end{proof}

\begin{corollary}\label{unramified}
If $\sE$  is unramified over $\sR$,  then
$\SK(\sE,\tau) \cong \SK(\sE_0,\tau|_{\sE_0})$.
\end{corollary}

\begin{proof}
Since $\sE$ is unramified over $\sR$, we have $\sT$ is unramified over
$\sR$, $Z(\sE_0) = \sT_0$,  and  $\Gamma_\sE = \Gamma_\sR$,
so we can choose all the $x_\gamma$'s to lie in $\sR$.
 The assertion thus follows immediately
 from Th.~\ref{msem}, as $P = \Sigma_{\tau|_{\sE_0}}(\sE_0)$
and $X \subseteq R_0^* \subseteq P$.
(Alternatively, more  directly, one can
observe that $\Sigma'_0=\Sigma'_{\tau|_{\sE_0}}(\sE_0)$ and
$\Sigma_0=\Sigma_{\tau|_{\sE_0}}(\sE_0)$ and so deduce the Corollary by
Lemma~\ref{unramfacts}(\ref{nine}).)
\end{proof}

\begin{corollary}\label{seses}
If $\sT$ is unramified over $\sR$ and $\sE$ has a maximal graded subfield
$\sM$ unramified over $\sT$ and another maximal graded subfield $\sL$ totally
ramified over $\sT$  with $\tau(\sL ) =\sL$, then $\sE$ is semiramified with
$\sE_0=\sM_0$ $($a field$)$ and $\Gamma_\sE=\Gamma_\sL$, and
\begin{equation}\label{semirsk}
\SK(\sE,\tau) \ \cong \ \big\{a \in \sE_0
\mid N_{\sE_0/\sT_0}(a)\in \sR_0\big\}
 \, {\big/} \,  \textstyle{\prod\limits_{h\in \Gal(\sE_0/\sT_0)}}
\sE_0^{*h\overline \tau}.
\end{equation}
\end{corollary}
\begin{proof}
Let $n=\ind(\sE)$. Since $[\sM_0:\sT_0]=n$ and
${|\Gamma_\sL:\Gamma_\sT|=n}$, it follows from the Fundamental Equality
~\eqref{fundeq}   for $\sE/\sT$,
$\sM/\sT$, and $\sL/\sT$
that
$[\sE_0:\sT_0]=[\Gamma_\sE:\Gamma_\sT]=n$,  $\sE_0=\sM_0$,
which is
is a field, and $\Gamma_\sE = \Gamma_\sL$.
Thus $\sE$ is semiramified, so for the $\dlambda$ of
\eqref{dlambda},  $\dlambda = 1$. Now, $\sL^\tau$ is a graded
subfield of
$\sL$ with $[\sL:\sL^\tau] = 2$.  Since $\sL_0 = \sT_0$ while $(\sL^\tau)_0
= (\sL_0)^\tau = \sR_0$, $\sL$~must be unramified over $\sL^\tau$;
hence, $\Gamma_\sE = \Gamma_\sL  = \Gamma_{\sL^\tau}$.
Therefore, one can choose all the
$x_\gamma$'s in Th.~\ref{msem} to lie in $\sL^\tau$. Then each
 $x_\gamma x_\delta x_{\gamma+\delta}^{-1} \in (\sL^\tau)_0^* =
\sR_0^*=
\sE_0^{*\ov \tau}$. Hence, $X \subseteq P =
 \prod_{h\in \Gal(\sE_0/\sT_0)}\sE_0^{*h\overline \tau}$,
so
the formula for
 $\SK(\sE,\tau)$ in Th.~\ref{msem} reduces to (\ref{semirsk}).
\end{proof}

\begin{remark} In a sequel to this paper~\cite{II}, the following will be
shown:  With the hypotheses  of Th.~\ref{msem}, suppose
 $\sE$ is  semiramified with a graded  maximal subfield
$\sL$ totally ramified over $\sT$ such that $\tau(\sL) = \sL$, and
suppose
$\Gal(\sE_0/\sT_0)$ is bicyclic,
say $\sE_0 = N\otimes _{T_0} N'$ with
$N$ and $N'$ each cyclic Galois over $\sT_0$. Then,
$$
{\SK(\sE,\tau) \ \cong \ \Br(\sE_0/\sT_0;\sR_0)\big/
\big(\Br(N/\sT_0;\sR_0)+\Br(N'/\sT_0;\sR_0)\big)},
$$
 where
$\Br(\sE_0/\sT_0)$ is the relative Brauer group
$\ker\big(\Br(\sT_0)\to \Br(\sE_0)\big)$ and
$\Br(\sE_0/\sT_0;\sR_0)$ is the kernel of
$\text{cor}_{\sE_0\to \sR_0}\colon \Br(\sE_0/\sT_0)\to \Br(\sE_0/\sR_0)$.
(Compare this with~\cite[Th.~5.6]{y}.)  A further formula will be given
assuming only that $\sE$ is semiramified over $\sT$ with
$\Gal(\sE_0/\sT_0)$  bicyclic.
\end{remark}

In his construction of division algebras $D$ with nontrivial $\SK$,
Platonov worked
originally in~\cite[\S4]{platonov} with a division algebra $D$
where $Z(D)$ is a
Laurent power series field; he gave an exact sequence relating
$\SK(D)$ with
$\SK(\overline D)$ and what he called the ``group of projective
conorms.''
Yanchevski\u\i ~gave in~\cite[4.11]{y} an analogous exact
sequence in the unitary case.
Their results and proofs are valid whenever $Z(D)$ has a
henselian discrete
(rank $1$) valuation. We show here that their results hold more generally
whenever $Z(D)$ has a henselian
valuation with $\Gamma_D/\Gamma_{Z(D)}$ cyclic. We work in
the equivalent graded setting where the arguments are more transparent.

As before, let $\sE$ be a graded division algebra finite dimensional
over its center $\sT$ with a unitary graded involution~$\tau$, and let
$\sR=\sT^\tau$. Assume that $\sT$ is unramified over $\sR$ and  that
$\Gamma_\sE/\Gamma_\sT$ is cyclic group. (This cyclicity holds, e.g.,
 whenever $\Gamma_\sT\cong \mathbb Z$.)
It follows that the surjective map $\Theta_\sE\colon\Gamma_\sE
\rightarrow \Gal(Z(\sE_0)/\sT_0)$ has kernel~$\Gamma_\sT$. (For,
by \cite[Prop.~2.1, (2.3), Remark 2.4(i)]{hwcor}$, \ker(\Theta_\sE)/\Gamma_\sT$ has a nondegenerate symplectic
pairing, and hence has even rank as a finite abelian group. But here
$\ker(\Theta_\sE)/\Gamma_\sT$
is a cyclic group.) Hence,  $\dlambda=1$
by Lemma~\ref{dlambdafacts}, so $\sE$ is inertially split.
Invoking Lemma~\ref{unramfacts}\eqref{seven},
choose any $s\in \sE^*$
with $\deg(s)+\Gamma_\sT$ a generator of $\Gamma_\sE/\Gamma_\sT$,
such that $\tau(s)=s$. Let $\sigma=\va_s\in \Aut_\sT(\sE)$; so
$\sigma\tau$ is another $\sT/\sR$-graded involution of $\sE$, and
$\tau\sigma=\sigma^{-1}\tau$ (see Remark~\ref{easyob}(\ref{three})).
By the choice of $s$, $\sigma|_{Z(\sE_0)}$ is a generator of the
cyclic group $\Gal(Z(\sE_0)/\sT_0)$. Note that
$\Gal(Z(\sE_0)/\sR_0)=\langle \sigma|_{Z(\sE_0)},\tau|_{Z(\sE_0)}\rangle$
is a dihedral group.
Recall our convention that $c^{\sigma\tau}$ means $\sigma(\tau(c))$.
Let
\begin{align*}
\mathcal S \ & = \  \big \{ (\beta,b) \in Z(\sE_0)^* \times \sE_0^* \mid
\beta^{\sigma-1}=\Nrd_{\sE_0}(b)\big \} ;\\
\mathcal N \ & = \  \pi_1(\mathcal S)
\textrm{ (projection into the first component) } \\\ & \quad= \
\big \{ \beta \in Z(\sE_0)^* \mid  \beta^{\sigma-1}=
\Nrd_{\sE_0}(b) \textrm{ for some $b \in \sE_0^*$}\big \};\\
\mathcal W \ & = \ \sT_0^* \cdot \Nrd_{\sE_0}(\sE_0^*)  \ \subseteq
 \ Z(\sE_0)^*;\\
\mathcal P \ & = \ \mathcal N /\mathcal W ,\textrm{ which is Platonov's
group of projective conorms for $\sE$ }\textrm{\cite[\S4]{platonov}}.\qquad
\qquad\qquad\quad
\end{align*}
\begin{align*}
\mathcal S_\tau \ & = \ \big \{(\alpha,a) \in Z(\sE_0)^*\times\sE_0^* \mid
\alpha^{\sigma-1}=\Nrd_{\sE_0}(a)^{1-\sigma\tau} \big \};\\
\mathcal N_\tau \ & = \ \pi_1(\mathcal S_\tau)
\textrm{ (projection into the first component) } \\
\ & \quad = \ \big \{ \alpha \in Z(\sE_0)^* \mid \alpha^{\sigma-1}=
\Nrd_{\sE_0}(a)^{1-\sigma\tau} \textrm{ for some } a\in \sE_0^* \big \};\\
\mathcal W_\tau\ & = \ \sT_0^*\cdot \Nrd_{\sE_0}(\Sigma_\tau(\sE_0)) \
\subseteq \  Z(\sE_0)^* ;\\
\mathcal PU_\tau \ & = \ \mathcal N_\tau /\mathcal W_\tau, \textrm{ which is
Yanchevski\u\i's group of unitary projective conorms for }
(\sE,\tau)\textrm{~\cite[4.11]{y}}.
\end{align*}


\begin{proposition}\label{cyclic}
If $\sT$ is unramified over $\sR$ and $\Gamma_\sE/\Gamma_\sT$ is cyclic,
then for any generator
$\sigma$ of the cyclic group $\Gal(Z(\sE_0)/\sT_0)$,
we have

\begin{enumerate}[\upshape(i)]
\item \label{rh1} $\SK(\sE) \ \cong \  \{ a\in E_0^* \mid
N_{Z(\sE_0)/\sT_0}(\Nrd_{\sE_0}(a))=1 \}\, \big / \, \big([\sE_0^*,\sE_0^*]\cdot
\{ c^{\sigma-1}\mid c\in \sE_0^* \}\big).$

\smallskip

\item \label{rh2} $\SK(\sE,\tau) \ \cong  \  \{ a\in E_0^* \mid N_{Z(\sE_0)/\sT_0}
(\Nrd_{\sE_0}(a)) \in \sR_0  \} \, \big / \, \big(
\Sigma_\tau(\sE_0)\cdot \Sigma_{\sigma\tau}(\sE_0)\big)$.

\smallskip

\item \label{rh3} The following sequence is exact:
\begin{equation}
\SK(\sE_0,\sigma\tau)  \longrightarrow \SK(\sE,\tau)
\overset {f}{\longrightarrow}  \mathcal PU_\tau \longrightarrow 1,
\end{equation}
where the map $f\colon\SK(\sE,\tau)\rightarrow \mathcal PU_\tau$ is
the composition of $a\, \Sigma_\tau(\sE_0)\cdot \Sigma_{\sigma\tau}(\sE_0)
\mapsto (\alpha,a) \in \mathcal S_\tau$ and
${(\alpha,a)\mapsto \alpha \mathcal W_\tau \in \mathcal PU_\tau}$.
\smallskip

\item \label{rh4} There is a commutative diagram with exact rows:
\begin{equation*}
\begin{split}
\xymatrix{
& \SK(\sE_0,\sigma\tau) \ar[r] \ar[d]_{{}_{\substack{a\\ \downarrow\\ a^{1-\sigma\tau}}}}
 & \SK(\sE,\tau) \ar[r]^-{f} \ar[d]_{{}_{\substack{a\\ \downarrow\\ a^{1-\sigma\tau}}}}
 & \mathcal PU_\tau \ar[r] \ar[d]_{{}_{\substack{\alpha\\ \downarrow\\ \alpha}}}&  1\\
1 \ar[r] & \SK(\sE_0) \ar[r] \ar[d]_{{}_{\substack{b\\ \downarrow\\ b}}}
 & \SK(\sE) \ar[r]^-{g} \ar[d]_{{}_{\substack{b\\ \downarrow\\ b}}}
 & \mathcal P \ar[r] \ar[d]_{{}_{\substack{\beta\\ \downarrow\\ \beta^{1+\tau}}}} & 1\\
& \SK(\sE_0,\sigma\tau) \ar[r] &  \SK(\sE,\tau) \ar[r] & \mathcal PU_\tau \ar[r] & 1}
\end{split}
\end{equation*}
where the
 map $g\colon\SK(\sE)\rightarrow \mathcal P$ is the composition of
$ \ b\, [\sE_0^*,\sE_0^*]\langle c^{\sigma-1}\rangle \mapsto (\beta,b)\in \mathcal S$
and ${(\beta,b)\mapsto \beta \mathcal W \in \mathcal P}$.
\vskip.06truein
\item \label{unrkl} If  $\sE_0$ is a field, then $\SK(\sE,\tau)=1$.
\end{enumerate}
\end{proposition}
\begin{proof}  \hfill

(i) This formula was given by Suslin~\cite[Prop.~1.7]{sus1} for a division algebra
over a field  with a complete discrete valuation. In order to prove it in the
graded setting we need  two exact sequences which were given in
 \cite[Th.~3.4]{hazwadsworth}:
$$
\Gamma_\sE/\Gamma_\sT  \wedge \Gamma_\sE  /\Gamma_\sT
\longrightarrow \sE^{(1)}\big /[\sE_0^*,\sE^*]
\longrightarrow \SK(\sE) \longrightarrow 1,
$$
$$
1 \longrightarrow \ker \widetilde N /[\sE_0^*,\sE^*] \longrightarrow
\sE^{(1)}\big /[\sE_0^*,\sE^*] \longrightarrow
\mu_\dlambda(\sT_0) \cap \widetilde N(\sE_0^*)\longrightarrow 1,
$$
where $\sE^{(1)} = \{a\in \sE^*\mid \Nrd_\sE(a) = 1\} \subseteq \sE_0$ and
$\widetilde N = N_{Z(\sE_0)/\sT_0}\circ \Nrd_{\sE_0}\colon
\sE_0^* \to \sT_0^*$. Since $\dlambda=1$ (see the paragraph prior to the
Proposition) and the wedge product of a cyclic group with itself is trivial,
these  exact sequences yield
$$
\SK(\sE) \ \cong \  \{ a\in \sE_0^* \mid
N_{Z(\sE_0)/\sT_0}(\Nrd_{\sE_0}(a))=1 \}\big /[\sE_0^*,\sE^*].
$$
We are left to show that $[\sE_0^*,\sE^*]=[\sE_0^*,\sE_0^*]\cdot
\{ c^{\sigma-1}\mid c\in \sE_0^* \}$. This follows from  the fact
that $\sE^*/\sT^* \sE_0^*\cong \Gamma_\sE/\Gamma_\sT$ is cyclic
together with  the
following observation, which is easily verified  using the standard commutator
identities: If $G$ is a group and $N$ is a normal subgroup of $G$ such that
$G/Z(G)N$ is a cyclic group generated by, say,  $xZ(G)N$,
then $[N,G]=[N,N][x,N]$ where $[x,N]=\{[x,n] \mid n \in N\}$.
(Here, take $G = \sE^*$, $N= \sE_0^*$,  and for $x$ take any $s\in \sE_\gamma^*$
for any $\gamma\in \Gamma_\sE$ such that $\Theta_\sE(\gamma+\Gamma_\sT)
= \sigma$.)

(ii) By Th.~\ref{msem}, taking into account that $\dlambda=1$ and
$\Gal(Z(\sE_0)/\sT_0)=\langle \sigma \rangle$, we have,
\begin{equation} \label{cycliccase}
\begin{split}
\SK(\sE,\tau)\ &  \cong  \ (\Sigma'_\tau \cap \sE_0) \big /
(\Sigma_\tau \cap \sE_0)   \\
\ & = \  \{ a \in \sE_0^*\mid
N_{Z(\sE_0)/\sT_0}\Nrd_{\sE_0}(a)\in \sR_0  \}\,
\big / \, \big(\Sigma_\tau(\sE_0)\cdot \Sigma_{\sigma\tau}(\sE_0)\cdot
\langle x_\gamma x_\delta x_{\gamma+\delta}^{-1}
\mid \gamma, \delta \in \Gamma_\sE \rangle\big),
\end{split}
\end{equation}
where for each $\gamma \in \Gamma_\sE$,
$x_\gamma$ is chosen in  $\sE^*_\gamma$
with $x_\gamma=\tau(x_\gamma)$ and $x_\gamma\ne 0$,
using Lemma~\ref{unramfacts}\eqref{seven}.
Let $\sL=\sR[s]$, where $s$ is chosen in $\sE^*$ with $\varphi(s)|_{\sE_0}
= \sigma$, which is possible as $\Theta_{\sE}\colon \Gamma_\sE \to
 \Gal(Z(\sE_0)/\sT_0)$
is surjective (see \eqref{surj}).  Moreover, $s$ can  be chosen  with
$\tau(s) = s$. Since $\ker(\Theta_\sE) = \Gamma_\sT = \Gamma_\sR$,
we have $\Gamma_\sE = \langle \deg(s) \rangle+ \Gamma_\sR$.  Thus,
$\sL$  is a graded
subfield of $\sE$ with $\Gamma_\sL=\Gamma_\sE$ and $\tau|_\sL=\id$.
For each $\gamma \in \Gamma_\sE$
we can choose $x_\gamma \in \sL_\gamma^*$; then for all
$\gamma,\delta \in \Gamma_\sE$, we have
$x_\gamma x_\delta x_{\gamma+\delta}^{-1} \in \sL_0^*
\subseteq \Sigma_\tau(\sE_0)$. Thus,  the
$\langle x_\gamma x_\delta x_{\gamma+\delta}^{-1}  \rangle$ term in
\eqref{cycliccase}  is redundant,
yielding the formula in (ii).

(iii) We first check that $f$ is well-defined: Take any $a \in \sE_0^*$
with $N_{Z(\sE_0)/\sT_0}(\Nrd_{\sE_0}(a))\in \sR_0^*$. Let
$c=\Nrd_{\sE_0}(a)$. Then, as $\sR_0=\sT_0^\tau$,
$1=N_{Z(\sE_0)/\sT_0}(c)^{1-\tau}=N_{Z(\sE_0)/\sT_0}(c)^{1-\sigma\tau}=
N_{Z(\sE_0)/\sT_0}(c^{1-\sigma\tau})$. By Hilbert~90, there is
$\alpha \in Z(\sE_0)^*$ with $\alpha^{\sigma-1}=c^{1-\sigma\tau}=
\Nrd_{\sE_0}(a)^{1-\sigma\tau}$. Hence $(\alpha,a)\in \mathcal S_\tau$, so
$\alpha \in \mathcal N_\tau$, and the choice of $\alpha$ is unique up to
$\sT_0^* \subseteq \mathcal W_\tau$.
Thus, the image of $a$ in $\mathcal  P U_\tau$ is independent of the choice of
$\alpha$.
Suppose further that $a=pq$ for some
$p\in \Sigma_\tau(\sE_0)$, $ q \in \Sigma_{\sigma\tau}(\sE_0)$, say,
$p=s_1\ldots s_k$ with each $s_i \in S_\tau(\sE_0)$.
Then,
\begin{equation}\label{skks}
\begin{split}
\Nrd_{\sE_0}(p)^\tau \ &= \ \Nrd_{\sE_0}(s_1)^\tau \ldots
\Nrd_{\sE_0}(s_k)^\tau  \  = \ \Nrd_{\sE_0}(s_1^\tau)\ldots
\Nrd_{\sE_0}(s_k^\tau)\\
& = \ \Nrd_{\sE_0}(s_1)\ldots \Nrd_{\sE_0}(s_k) \ = \ \Nrd_{\sE_0}(p);
\end{split}
\end{equation}
likewise, $\Nrd_{\sE_0}(q)^{\sigma\tau}=\Nrd_{\sE_0}(q)$. So,
$$
\alpha^{\sigma-1} \, = \ \Nrd_{\sE_0}(pq)^{1-\sigma\tau} \,
= \ \Nrd_{\sE_0}(p)^{1-\sigma\tau}\Nrd_{\sE_0}(q)^{1-\sigma\tau} \,
= \ \Nrd_{\sE_0}(p)^{1-\sigma}.
$$
Hence, $\big (\alpha\Nrd_{\sE_0}(p)\big )^{\sigma-1}=1$, showing that
$\alpha \Nrd_{\sE_0}(p) \in \sT_0$; Thus $\alpha \in \mathcal W_\tau$.
This proves that $f$ is well-defined.

For the subjectivity of $f$, take any $\alpha \in \mathcal N_\tau$.
Then, there is $a\in \sE_0^*$ with
$\alpha^{\sigma-1}=\Nrd_{\sE_0}(a)^{1-\sigma\tau}$. So,
$N_{Z(\sE_0)/\sT_0}(\Nrd_{\sE_0}(a))^{1-\sigma\tau}
=N_{Z(\sE_0)/\sT_0}(\alpha^{\sigma-1})=1$, which shows that
$N_{Z(\sE_0)/\sT_0}(\Nrd_{\sE_0}(a))\in \sT_0^{\sigma\tau}
=\sT_0^\tau=\sR_0$, and hence
$a \in \Sigma_\tau'(\sE)\cap \sE_0^*$. Since
$f\big(a\Sigma_\tau(\sE_0)\Sigma_{\sigma\tau}(\sE_0)\big)
=\alpha\mathcal W_\tau$, $f$ is surjective.

Finally, we determine $\ker(f)$: The image of $\SK(\sE_0,\sigma\tau)$ in
$\SK(\sE,\tau)$ is $\Sigma_{\sigma\tau}'(\sE_0)\Sigma_\tau(\sE_0)\big /
\Sigma_{\sigma\tau}(\sE_0)\Sigma_\tau(\sE_0)$. An element in this image
is represented by some
$a \in \Sigma_{\sigma\tau}'(\sE_0)$. For such an $a$,
$\Nrd_{\sE_0}(a)^{1-\sigma\tau}=1$. Then $(1,a) \in \mathcal S_\tau$, so
that $f$ maps the image of $a$ to $1$ in $\mathcal PU_\tau$. Conversely,
suppose
$a\Sigma_\tau(\sE_0)\Sigma_{\sigma\tau}(\sE_0) \in \ker(f)$. That is,
$\Nrd_{\sE_0}(a)^{1-\sigma\tau}=\alpha^{\sigma-1}$, where
$\alpha \in \mathcal W_\tau$,  so
$\alpha=c\Nrd_{\sE_0}(d)$ with $c\in \sT_0^*$ and $d\in \Sigma_\tau(\sE_0)$.
So, $\Nrd_{\sE_0}(d)=\Nrd_{\sE_0}(d)^{\tau}$ by the argument
of~(\ref{skks}) above, and hence
$$
\Nrd_{\sE_0}(a)^{1-\sigma\tau} \, = \ \alpha^{\sigma-1}
 \, = \ \big (c\Nrd_{\sE_0}(d)\big )^{\sigma-1}
 \, = \ \Nrd_{\sE_0}(d)^{\sigma-1} \, = \ \Nrd_{\sE_0}(d)^{\sigma\tau-1}.
$$
Thus, $\Nrd_{\sE_0}(ad)^{1-\sigma\tau}=1$, i.e.,
$ad \in \Sigma_{\sigma\tau}'(\sE_0)$. Hence,
$a=(ad)d^{-1} \in \Sigma_{\sigma\tau}'(\sE_0)\Sigma_\tau(\sE_0)$. This
shows that $\ker(f)$ coincides with the image of $\SK(\sE_0,\sigma\tau)$
in $\SK(\sE,\tau)$, completing the proof of exactness of the  sequence.

(iv) Exactness of the middle row is proved by an analogous but easier
argument to that for~(\ref{rh3}).
Commutativity of the left  rectangles of the diagram is evident.
Commutativity of
the top right rectangle  is clear from the definitions.
Commutativity of the bottom right rectangle is
easy to check using the identity
\begin{equation}\label{sigmatau}
(1-\sigma\tau)\circ(\sigma-1) \ = \ (\sigma-1)\circ(1+\tau),
\end{equation}
which follows from $(\sigma\tau)^2
 = \id$. Note that for each column of the diagram,
the composition of the two maps is the squaring map.

(v) For this part, the proof follows closely Yanchevski\u\i's proof in
\cite[4.13]{y}. (But our notational convention for products of functions is
$fg = f\circ g$, whereas his appears to be $fg = g\circ f$.)
Suppose $\sE_0$ is a field. For simplicity we denote
$\overline \tau=\tau|_{\sE_0}$ by $\tau$. Take $a\in
\Sigma_\tau'(\sE)\cap \sE_0$. So,
$N_{\sE_0/\sT_0}(a^{1-\tau})=1$. We will
show that $a \in \sE_0^\tau \sE_0^{\sigma\tau}$.
It then follows by (\ref{rh2}) above
 that $\SK(\sE,\tau)=1$. But since $\sE_0$ is cyclic over
$\sT_0$, by  Hilbert~90 there is a $b \in \sE_0^*$ such that
$a^{\tau-1}=b^{\sigma-1}$ where $\langle \sigma \rangle  =\Gal(\sE_0/\sT_0)$.
So, ${1=a^{(\tau+1)(\tau-1)}=b^{(\tau+1)(\sigma-1)}}$.
Analogously to \eqref{sigmatau}, we have
 $(\tau+1)(\sigma-1)=(\sigma-1)(1-\tau\sigma)$.
So $b^{(\sigma-1)(1-\tau\sigma)}=1$.
Setting $c=b^{(1-\tau\sigma)}$, we have $c^{\sigma-1}=1$, so  $c\in \sT_0$.
But, $N_{\sT_0/\sR_0}(c)=c^{1+\tau\sigma}=
b^{(1+\tau\sigma)(1-\tau\sigma)}=1$.  By Hilbert 90  we have
$c=d^{\tau\sigma-1}$ for some $d\in \sT_0^*$. Let $t=bd\in  \sE_0^*$. Then,
$t^{1-\tau\sigma}=b^{(1-\tau\sigma)}d^{(1-\tau\sigma)}=
d^{(\tau\sigma-1)}d^{(1-\tau\sigma)}=1$, i.e.,  $t \in \sE_0^{\tau\sigma}$.
So, $\sigma(t) = \tau(t) \in \sE_0^{\sigma\tau}$.
Thus,  $a^{\tau-1}=b^{\sigma-1}=(t/d)^{\sigma-1}=t^{\sigma-1}=
t^{\tau-1}$, as $d\in \sT_0$. This shows that
$(a\tau(t))^{\tau-1}=1$, i.e., $a\tau(t) \in \sE^\tau$; hence
$a=(a\tau(t))\tau(t)^{-1}\in  \sE_0^\tau
\sE_0^{\sigma\tau}$.
\end{proof}

\section{Totally ramified algebras}

For a graded division algebra $\sE$ totally ramified over its center
$\sT$ with a unitary graded involution $\tau$, two possible cases can
arise: either $\sT$ is totally ramified over $\sR=\sT^\tau$, or $\sT$ is
unramified over $\sR$. In the first case, we showed in Prop.~\ref{total}
that $\SK(\sE,\tau)$ is trivial. We now obtain an
easily computable explicit formula
for $\SK(\sE,\tau)$ in the second case.
For a field $K$ and for $n\in \nn$, we write $\mu_n$ for the group of
all $n$-th roots of unity in an algebraic closure of $K$.  Then set
$\mu_n(K) = \mu_n \cap K^*$.

\begin{theorem} \label{sktotal} If $\sE$ is totally ramified over $\sT$ of
index $n$ and $\sT$ is unramified over $\sR$, then
\begin{align}
\SK(\sE,\tau)\ & \cong  \  \big \{a\in \sT_0^*\mid a^n\in \sR_0^*\} \,\big /
 \, \{a\in \sT_0^*\mid a^e\in \sR_0^* \} \label{genesa}
\\
&\cong  \  \big \{\omega \in \mu_n(\sT_0) \mid \tau(\omega)
=\omega^{-1}\big \}\big /\mu_e, \label{genesa1}
\end{align}
where $e$ is the exponent of $\Gamma_\sE/\Gamma_\sT$. In particular,
\begin{enumerate}[\upshape(i)]
\item \label{hjh1} The restriction of the map
$K_1(\sE,\tau)\rightarrow K_1(\sE)$ given by  $a\Sigma_\tau \mapsto a^{1-\tau}\sE'$, induces an  injective map
$$
\alpha\colon \SK(\sE,\tau)\longrightarrow \SK(\sE)\ \cong \
\mu_n(\sT_0)/\mu_e.
$$
\item \label{hhh} If the exponent $e$ of $\sE$ is odd, then $\alpha$ is an
isomorphism.

\smallskip

\item  If  $e>2$ then $\sT_0=\sR_0(\mu_e)$, and $\tau$ acts on $\mu_e$ by
$\omega\mapsto \omega^{-1}$.

\end{enumerate}
\end{theorem}
\begin{proof}
Since $\sT$ is unramified over $\sR$ and $\sE_0=\sT_0$, the formulas of
Th.~\ref{msem} for $\SK(\sE,\tau)$ reduce to $\dlambda = n$ and
\begin{equation}\label{gendesb}
\SK(\sE,\tau) \ \cong \  \{a\in \sT_0^*\mid a^n\in \sR_0^* \} \,
\big / \, \big(\sR_0^*  \,  \langle x_\gamma x_\delta x_{\gamma+\delta}^{-1}
\mid \gamma, \delta \in \Gamma_\sE \rangle\big),
\end{equation}
where each $x_\gamma\in \sE_\gamma^*$ with $\tau(x_\gamma) = x_\gamma$.
Recall that as $\sE/\sT$ is totally ramified,
the canonical pairing $\sE^*\times \sE^* \rightarrow \mu_e(\sT_0)$
given by  $(s,t)\mapsto [s,t]$ is surjective
(\cite[Prop.~2.1]{hwcor}), and   $\mu_e(\sT_0) = \mu_e$,
i.e., $\sT_0$ contains all $e$-th roots of unity.
Since each $\sE_\gamma = \sT_0 x_\gamma$ with $\sT_0$ central,
it follows  that
$\{[x_\delta,x_\gamma] \mid \gamma, \delta  \in \Gamma_\sE\} = \mu_e$.
Now consider $c=x_\gamma x_\delta x_{\gamma+\delta}^{-1}$ for any $\gamma,
\delta\in\Gamma_\sE$. Then,  $\tau(c)=x_{\gamma+\delta}^{-1}x_\delta x_\gamma$.
Note that $x_\delta x_\gamma$ and $x_{\gamma+\delta}$
each lie in $\sE_{\gamma+\delta} = \sT_0x_{\gamma+\delta}$,
so they commute.  Hence,
\begin{equation}\label{c}
\tau(c) c^{-1} \, = \ x_{\gamma+\delta}^{-1} (x_\delta x_\gamma)
x_{\gamma+\delta}x_\delta^{-1}x_\gamma^{-1}
 \ = \  [x_\delta,x_\gamma].
\end{equation}
  Since  $[x_\delta,x_\gamma] \in \mu_e$,
this shows that
$c \in \big \{a\in \sT_0^*\mid a^e\in \sR_0^*\big \}$. For the reverse
inclusion, take any
$d$~in~$\sT_0^*$ such that $d^e\in \sR_0^*$.
So $\tau(d)d^{-1} \in \mu_e$. Thus,
$\tau(d) d^{-1}= [x_\delta,x_\gamma]$, for some $\gamma, \delta\in \Gamma_\sE$.
 Taking ${c=x_\gamma x_\delta x_{\gamma+\delta}^{-1}}$, we have
$\tau(d) d^{-1}=\tau(c) c^{-1}$
by \eqref{c},  which implies that $dc^{-1}$ is
$\tau$-stable, so lies in $\sR_0^*$;
 thus, ${d \in \sR_0^*  \,  \langle x_\gamma x_\delta x_{\gamma+\delta}^{-1}
\mid \gamma, \delta \in \Gamma_\sE \rangle}$. Therefore,
$\sR_0^*  \,  \langle x_\gamma x_\delta x_{\gamma+\delta}^{-1} \mid
\gamma, \delta \in \Gamma_\sE \rangle
=\{a\in \sT_0^*\mid a^e\in \sR_0^* \}$. Inserting this in ~(\ref{gendesb}) we
obtain (\ref{genesa}).

(i) Consider the well-defined map
$\alpha\colon\SK(\sE,\tau)\rightarrow \SK(\sE)$ given by
$a\Sigma_\tau\mapsto a^{1-\tau}\sE'$ (see  diagram~\eqref{goodd} for the
non-graded version).
By~\cite[Cor.~3.6(ii)]{hazwadsworth},
$\SK(\sE)\cong \mu_n(\sT_0)/\mu_e$. Taking into account
formula ~\eqref{genesa} for $\SK(\sE,\tau)$,
it is easy to see that $\alpha$ is  injective.

We now verify that
\begin{equation}\label{imalpha}
\im(\alpha) \ = \ \big \{\omega \in \mu_n(\sT_0) \mid
\tau(\omega)=\omega^{-1}\big \} \, \big /\mu_e,
\end{equation}
and thus obtain ~(\ref{genesa1}).
Indeed, since $\mu_e =  \{[x_\delta,x_\gamma] \mid \gamma, \delta
\in \Gamma_\sE\}$,  by setting
$c=x_\gamma x_\delta x_{\gamma+\delta}^{-1}$ we have
${[x_\delta,x_\gamma]=\tau(c) c^{-1}}$ by \eqref{c}. This shows that $\mu_e \subseteq \big \{\omega \in \mu_n(\sT_0) \mid \tau(\omega)=\omega^{-1}\big \}$. Now for
any $\omega \in \mu_n(\sT_0)$ with $\tau(\omega)=\omega^{-1}$,
we have $N_{\sT_0/\sR_0}(\omega)=1$, so  Hilbert~90 guarantees that
$\omega=c^{1-\tau}$ for some $c\in \sT_0^*$. Then,
$(c^n)^{1-\tau} = \omega^n = 1$, so $c^n\in \sR_0^*$.
Thus,
$c\in \Sigma_\tau'$, and clearly  $\alpha(c\Sigma_\tau)=\omega\mu_e$.
This shows $\supseteq$ in \eqref{imalpha}; the reverse inclusion is clear
from  the definition of $\alpha$.

(ii)
Suppose $e$ is odd.  Let $m = |\mu_n(\sT_0)|$.  So, $\mu_n(\sT_0) = \mu_m$,
with $m \,|\, n$.  Also, $e\,|\, m$, as $\mu_e\subseteq \sT_0$.  Since $e$ and $n$
have the same prime factors, this is also true for $e$ and $m$.  Recall that
$\Aut(\mu_m)\cong (\zz/m\zz)^*$,  the multiplicative group of units of the ring
$\zz/m\zz$; so, $|\Aut(\mu_m)| = \varphi(m)$, where $\varphi$ is Euler's
$\varphi$-function. Since $e\,|\, m$ and $e$ and $m$ have the same prime
factors (all odd), the canonical map $\psi\colon \Aut(\mu_m) \to
\Aut(\mu_e)$ given by restriction is surjective with kernel of order
$\varphi(m)/\varphi(e) = m/e$, which is odd.  Therefore,
$\psi$ induces an isomorphism on the $2$-torsion subgroups,
$_2\!\Aut(\mu_m) \cong \ _2\!\Aut(\mu_e)$.  Now, $\tau|_{\mu_m}
\in \, _2\!\Aut(\mu_m)$ and we saw for~(i) that $\tau|_{\mu_e}$
is the inverse map $\omega \mapsto \omega^{-1}$.  The inverse map
on $\mu_m$ also lies in $_2\!\Aut(\mu_m)$ and has the same restriction
to $\mu_e$ as $\tau$.  Hence, $\tau|_{\mu_m}$ must be the inverse map.
That is, $\{\omega \in \mu_n(\sT_0) \mid
\tau(\omega)=\omega^{-1} \} = \mu_n(\sT_0)$.  Therefore,
\eqref{imalpha} above shows that $\im(\alpha) = \mu_n(\sT_0)/\mu_e$,
which we noted above is isomorphic to $\SK(\sE)$.

(iii) We saw in the proof of part (i) that
$\tau$ acts on $\mu_e$ by the inverse map.  So,
if $e>2$, then $\mu_e \not \subseteq \sR_0$.
Since $[\sT_0:\sR_0]=2$, it then follows that
$\sT_0=\sR_0(\mu_e)$.
\end{proof}

\begin{remark}
The isomorphism $\SK(\sE,\tau)\cong \SK(\sE)$ of part (\ref{hhh}) of the
above theorem can be obtained under the milder condition that
$\sE_0=\sT_0\sE'$ provided that the exponent of $\sE$ is a prime power. The
proof is similar.
\end{remark}

\begin{example} \label{toex}
Let $r_1, \ldots, r_m$ be integers with each
$r_i \ge 2$.  Let $e = \text{lcm}(r_1, \ldots, r_m)$, and let
$n = r_1\ldots r_m$.  Let $C$ be any field such that $\mu_e
\subseteq C$ and $C$ has an automorphism $\theta$
of order $2$ such that
$\theta(\omega) = \omega^{-1}$ for all $\omega\in \mu_e$.
Let $R$ be the fixed field $C^\theta$.
Let  $x_1,\dots, x_{2m}$ be $2m$ independent indeterminates, and let
$K$ be the  iterated Laurent power series field $C((x_1))\dots ((x_{2m}))$.
This $K$~is equipped with its standard valuation
$v\colon K^* \rightarrow \mathbb Z^{2m}$ where $\mathbb Z^{2m}$ is given
the right-to-left lexicographical ordering. With this valuation $K$ is
henselian (see~\cite[p.~397]{wadval}). Consider the tensor product of
symbol algebras
\begin{equation*}
D\ = \ \Big (\frac{x_{1},x_{2}}{K}\Big )_{\omega_1} \otimes_K
\ldots \otimes_K
\Big (\frac{x_{2m-1},x_{2m}}{K}\Big )_{\omega_m},
\end{equation*}
where  for $1\leq i \leq m$, $\omega_i$ is a primitive $r_i$-th root of unity
in $C$.
Using the valuation theory developed for division algebras, it is known that
$D$ is a division algebra, the valuation $v$ extends to $D$, and $D$ is
totally ramified over $K$ (see \cite[Ex.~4.4(ii)]{wadval} and
\cite[Ex.~3.6]{tw}) with
$$
\Gamma_D/\Gamma_K \ \cong \
\textstyle{\prod\limits_{i=1}^m} (\mathbb Z/r_i \mathbb Z)\times
(\mathbb Z/r_i\mathbb Z),
$$
and $\ov D = \ov K\cong  C$.
Extend $\theta$ to an automorphism  $\theta'$ of order $2$ on $K$ in the obvious
way, i.e., acting by~$\theta$  on the coefficients of a Laurent series, and with
$\theta'(x_i) = x_i$ for $1\le i\le 2m$.
On each of the symbol algebras
$\Big (\frac{x_{2i-1},x_{2i}}{K}\Big )_{\omega_i}$ with its generators $\ba_i$
and $\bb_i$ such that $\ba_i^{r_i}=x_{2i-1}$, $\bb_i^{r_i}=x_{2i}$, and
${\ba_i\bb_i=\omega_i\bb_i\ba_i}$, define an involution  $\tau_i$ as follows:
$\tau_i(c \, \ba_i^k \bb_i^l)= \theta'(c)  \, \bb_i^l \ba_i^k$, where $c\in K$ and $0\leq l,k < r_i$.
 Clearly
${K^{\tau_i}= K^{\theta'} = R((x_1))\dots ((x_{2m}))}$, and therefore $\tau_i$~is a
unitary involution.  Since the $\tau_i$ agree
on $K$ for $1\le i\le m$, they yield   a unitary involution
$\tau=\otimes_{i=1}^m\tau_i$
on $D$. Now by Th.~\ref{involthm2},
$\SK(D,\tau)\cong\SK(\gr(D),\widetilde \tau)$. Since
$D$ is totally ramified over $K$, which is unramified over $K^{\tau}$,
we have correspondingly that $\gr(D)$ is totally ramified over~$\gr(K)$,
which is unramified over $\gr(K)^{\widetilde \tau}$.
Also, $\gr(K)_0 \cong \ov K \cong  C$.
We have ${\exp(\gr(D)) = \exp(D) = \exp(\Gamma_D/\Gamma_K) =
\text{lcm}(r_1, \ldots, r_m) = e}$ and
$ {\ind}(\gr(D)) = \ind(D) = r_1\ldots r_m = n$.
By Th.~\ref{sktotal},
$$
\SK(D,\tau) \ \cong \ \SK(\gr(D),\widetilde \tau)\ \cong \
\{ \omega\in \mu_n(C) \mid \theta(\omega) = \omega^{-1}\}\big/ \mu_e,
$$
while by \cite[Th.~4.8, Cor.~3.6(ii)]{hazwadsworth},
$$
\SK(D) \ \cong \ \SK(\gr(D))\ \cong \mu_n(C)/\mu_e .
$$
Here are some more specific examples:

(i) Let $C = \mathbb C$, the complex numbers, and let
$\theta$ be complex conjugation, which maps every root of
unity to its inverse. So, $R = C^\theta = \mathbb R$.
Then, $\SK(D, \tau) \cong \SK(D) \cong \mu_n/\mu_e \cong\zz
/(n/e)\zz$.

(ii) Let $r_1 = r_2 = 4$, so $e = 4$ and $n = 16$.  Let
$\omega_{16}$ be a primitive sixteenth root of unity in $\cc$,
and let
${C = \qq(\omega_{16})}$, the sixteenth
cyclotomic extension of $\qq$.  Recall that ${\Gal(C/\qq)
 \cong \Aut(\mu_{16}) \cong (\zz/4\zz) \times (\zz/2\zz)}$,
Let $\theta\colon C \to C$ be the automorphism which
maps  $\omega_{16} \mapsto (\omega_{16})^7$. Then,
$\theta^2 = \id_C$, as $7^2 \equiv 1 \ (\text{mod }16)$,
and $\{\omega \in  \mu_{16}\mid \theta(\omega) = \omega^{-1}\}
 = \mu_8$.  Thus, $\SK(D, \tau) \cong \mu_8/\mu_4 \cong
\zz/2\zz$, while $\SK(D) \cong \mu_{16}/\mu_4 \cong \zz/4\zz$.
So, here the injection  $\SK(D,\tau) \to \SK(D)$ is not surjective.

(iii) Let $r_1 = \ldots = r_m = 2$, so $e = 2$ and $n = 2^m$.
Here, $C$ could be any quadratic extension of any field~
$R$ with $\chr(R) \ne 2$. Take $\theta$ to be the unique
nonidentity $R$-automorphism of $C$.  The resulting
$D$ is a tensor product of
$m$ quaternion algebras  over $C((x_1))\ldots((x_{2m}))$,
and ${\SK(D,\tau) \cong
\{\omega\in \mu_{2^m}(C)\mid \theta(\omega) = \omega^{-1}\}\big/
\mu_2}$, while $\SK(D) \cong \mu_{2^m}(C)/\mu_2$.
\end{example}

Ex.~\ref{toex} gives an indication how to use the graded approach to
recover results in the literature on
the unitary $\SK$ in a unified
manner and to extend them from division algebras with discrete valued groups to
 arbitrary valued groups.  While $\SK(D)$ has long been known for the $D$
of Ex.~\ref{toex},
the formula for $\SK(D,\tau)$ is new.

Here is a more complete statement of what the results in the preceding sections
 yield for
$\SK(D, \tau)$ for valued division algebras $D$.

\begin{theorem}\label{appl}
Let $(D,v)$ be a tame valued division algebra over a  field $K$
with $v|_K$ henselian,
with a unitary involution $\tau$; let $F=K^\tau$, and suppose $v|_F$ is
henselian  and that $K$ is tamely
ramified over $F$. Let $\ov \tau$ be the involution on $\ov D$
induced by $\tau$.
Then,
\begin{enumerate}
  \item[$(1)$] Suppose $K$ is unramified over $F$.
\begin{enumerate}[\upshape(i)]
  \item \label{apun} If $D$ is unramified over $K$, then
$\SK(D,\tau)\cong \SK(\overline D,\overline \tau)$.
\smallskip
  \item \label{tolast} If $D$ is totally ramified over $K$, let $e = \exp(D)$
and $n = \ind(D)$; then,
  $$
\SK(D,\tau) \ \cong  \  \{\omega \in \mu_n(\overline K) \mid \tau(\omega) =
\omega^{-1}\}\big /\mu_e,
$$
while $\SK(D) \cong \mu_n(\ov K)/\mu_e$.
\smallskip
  \item \label{tirl} If $D$ has a maximal graded subfield
$M$ unramified over $K$ and another maximal graded subfield $L$ totally
ramified over $K$, with $\tau(L ) =L$, then $D$ is semiramified  and
\begin{equation*}\label{semirsk3}
\SK(D,\tau) \ = \ \big\{a \in \overline D^* \mid N_{\overline D/\overline K}(a)\in \overline F\big\}
 \, {\big/} \,  \textstyle{\prod\limits_{h\in \Gal(\overline D/\overline K)}}
\overline F^{*h\overline \tau}.
\end{equation*}
  \item\label{gamcyclic} Suppose $\Gamma_D/\Gamma_K$ is cyclic.  Let
$\sigma$ be a generator of $\Gal(Z(\ov D)/\ov K)$.  Then,
$$
\SK(D, \tau) \ \cong \ \{ a\in \ov D^*\mid N_{Z(\ov D)/\ov K}(\Nrd_{\ov D}(a))
\in \ov F\} \, \big/ \, \big(\Sigma_{\ov \tau}(\ov D) \cdot
\Sigma_{\sigma\ov \tau}(\ov D)\big).
$$
  \item \label{apsem} If $D$ is  inertially split, $\overline D$ is a field and
$\Gal(\overline D/\overline K)$ is cyclic, then $\SK(D,\tau)=1$.
\end{enumerate}
\medskip
  \item[$(2)$] \label{apto} If $K$ is totally ramified over $F$, then
$\SK(D,\tau)=1$.
\end{enumerate}
\end{theorem}

\begin{proof}
Let $\gr(D)$ be the associated graded division algebra  of $D$.
The tameness assumptions assure that $\gr(K)$ is the center of
$\gr(D)$ with $[\gr(D):\gr(K)] = [D:K]$ and that the graded involution
$\widetilde\tau$ on $\gr(D)$ induced by $\tau$ is unitary with
$\gr(K)^{\widetilde \tau} = \gr(K^\tau)$. In each case of Th.~\ref{appl},
the conditions on $D$ yield analogous conditions on $\gr(D)$.
Since by Th.~\ref{involthm2}, $\SK(D,\tau)\cong \SK(\gr(D),\widetilde \tau)$,
(2) and (1)(v) follow immediately  from Prop.~\ref{completely}
and Prop.~\ref{cyclic}(\ref{unrkl}), respectively. Part (1)(i),
also follows from
Th.~\ref{involthm2}, and Cor.~\ref{unramified} as follows:
$$
\SK(D,\tau)  \ \cong \ \SK(\gr(D),\widetilde{\tau})
  \ \cong \ \SK(\gr(D)_0,\tau|_{\gr(D)_0}) \ = \
 \SK(\overline{D},\overline{\tau}).
$$
Parts (1)(ii), (1)(iii),  and (1)(iv) follow similarly using
Th.~\ref{sktotal}, Cor.~\ref{seses}, and Prop.~\ref{cyclic}\eqref{rh2}
respectively.
\end{proof}

In the special case  that the henselian valuation
on $K$ is discrete (rank $1$),
Th.~\ref{appl}~(1)(i), (iii), (iv), (v)  and (2) were obtained by
Yanchevski\u\i~ \cite{y}.  In this discrete case, the assumption
that $v$ on $K$ is henselian already implies that $v|_F$ is henselian
(see Remark~\ref{shensel}).

\end{document}